\documentclass{article}

\DeclareUnicodeCharacter{0301}{\'}
\usepackage{graphicx} %
\usepackage[]{biblatex}
\usepackage{setspace} %
\usepackage{geometry}
\usepackage{hyperref}
\usepackage{amsfonts}

\usepackage{tikz}
\usepackage{todonotes}
\usepackage{trackchanges}
\usepackage{algorithm}
\usepackage{algpseudocode}
\usepackage[normalem]{ulem} %
\usepackage{cancel} %
\usepackage{amsmath}

\usepackage{amssymb}
\usepackage{amsthm}

\usepackage{pict2e}
\usepackage{keyval}
\usepackage{calc}
\usepackage{fp}
\usepackage{diagbox}
\usepackage{booktabs}
\usepackage{tabularx,colortbl, makecell, caption, multirow}
\usepackage{boldline}
\usepackage[normalem]{ulem}

\newtheorem{proposition}{Proposition}
\newtheorem{theorem}{Theorem}

\newtheorem{observation}{Observation}

\newcommand{\R}{\mathbb{R}}

\def\sgn{\mbox{\rm sgn}}

\newtheorem{remark}[theorem]{Remark}

\usepackage{rotating}
\usepackage{graphicx}
\usepackage{bm}
\usepackage{xcolor}
\usepackage{multirow}
\usepackage{nicematrix}

\usepackage{graphicx}
\graphicspath{{images/}}
\usepackage{float}
\usepackage{longtable}
\usepackage{adjustbox}

\newcommand{\ds}{\displaystyle}

\usepackage{colortbl}

\graphicspath{{./figures/}}

\usepackage[titletoc,toc,title]{appendix}
\usepackage{comment}
\usepackage{graphics}
\usepackage{graphicx}
\usepackage{adjustbox}
\usepackage{tikz}
\usepackage{caption}
\usepackage{psfrag}
\usepackage{soul}
\usepackage{mwe}
\usepackage{pgfplots}
\pgfplotsset{compat=1.17} %

\graphicspath{{./figures/}}

\usepackage{comment}
\usepackage{graphics}
\usepackage{graphicx}
\usepackage{adjustbox}
\usepackage{tikz}
\usepgfplotslibrary{dateplot}
\usepgfplotslibrary{groupplots}
\usepackage{caption}
\usepackage{subcaption}
\usepackage{psfrag}
\usepackage{soul}
\usepackage{mwe}
\usepgfplotslibrary{fillbetween}

\newcommand{\fnote}{\textcolor{revblue}}
\usepackage{soul}
\definecolor{revblue}{RGB}{0,0,0}   %
\newcommand{\lnote}[1]{\textcolor{revblue}{#1}}
\newcommand{\lpnote}[1]{\textcolor{revblue}{#1}}   %

\definecolor{revgreen}{rgb}{0,0.45,0.15}   %
\newcommand{\gnote}[1]{\textcolor{revblue}{#1}}   %
\definecolor{yforange}{rgb}{0.85,0.45,0}

\usepackage{xspace}   %
\newcommand{\AOsub}{\fnote{\texttt{ACSO}\ensuremath{^{*}}}\xspace}   %
\newcommand{\AObase}{\fnote{\texttt{ACSO}}\xspace}                    %

\usepackage{authblk} %

\usepackage{color-edits}
\addauthor{yf}{purple}
\addauthor{fd}{blue}
\addauthor{lp}{violet}

\title{\textcolor{revblue}{ Feature Selection in Nonlinear SVMs via Local Search and Submodular
Optimization}}

\date{}

\author[1]{Federico D'Onofrio}
\author[2]{Yuri Faenza}
\author[1]{Laura Palagi}

\affil[1]{DIAG, Sapienza University of Rome}
\affil[2]{IEOR Department, Columbia University}

\begin{document}

\maketitle

\begin{abstract}

Embedded feature selection is a classical approach to interpretable machine learning, aiming to identify the most informative variables while simultaneously training the predictive model. We address the problem of enforcing an exact limit on the number of features used in nonlinear Support Vector Machines (SVMs) by means of a hard cardinality constraint. Whereas existing hard-constraint formulations apply only to the primal formulation of linear SVMs, we embed the cardinality constraint directly into the dual formulation of a nonlinear SVM. This ensures strict control over the number of selected features while still leveraging kernelization. We model the problem as a mixed-integer nonlinear program (MINLP). Our first contribution is a local-search metaheuristic applicable to any nonlinear kernel, which exploits the decomposable structure of the formulation to efficiently explore the binary decision space. Our second and main contribution is a decomposition framework that alternates optimization over continuous and binary variables. For polynomial kernels, we reformulate the binary subproblem as the maximization of a monotone submodular function under a cardinality constraint, enabling the integration of scalable greedy algorithms into the alternating scheme. Computational experiments on benchmark datasets show that the proposed algorithms significantly outperform existing MINLP solvers and standard heuristics, delivering higher-quality solutions within practical runtimes.

\end{abstract}

\section{Introduction}
Given the spread of  {opaque} Machine Learning (ML) models, the pursuit of inherently interpretable models has transitioned from a desirable attribute to a critical necessity \cite{rudin2022interpretable}. The lack of explanations for decisions made by complex ML models remains a key obstacle to their adoption in real‐world decision‐making scenarios \cite{goodman2017european, rudin2019stop}. Interpretable ML aims to develop models that offer insights into their decision‐making processes \cite{kamath2021explainable}. However, constructing such interpretable models often requires addressing complex optimization problems, which may demand significant computational resources, especially when dealing with large datasets. One effective strategy to enhance interpretability is \emph{Feature Selection} (FS), which consists of identifying and retaining only the most informative variables. {By enforcing sparsity, i.e., limiting the number of input features on which the model relies, we simplify the resulting decision rule, making it far more transparent and easier to interpret. } FS algorithms have gained significant traction in recent years because they not only enhance interpretability by reducing model complexity (\cite{blanquero2019variable,munoz2020informative}), but also clean the data and reduce noise (\cite{bolon2015feature,byeon2008simultaneously}), and improve predictive performance by mitigating overfitting (\cite{blanquero2020selection,lee2015mutual}). {By highlighting the most informative variables, FS yields models that end users  can more readily trust and understand.}
A comprehensive overview of various feature selection methods, including a discussion on their stability, is provided in \cite{chandrashekar2014survey}, while specific techniques have been developed for both regression \cite{andersen2010variable} and classification settings \cite{bertolazzi2016integer,tang2014feature}.

In this paper, we focus on the feature selection problem in \emph{Support Vector Machines} (SVMs). %
Among powerful ML tools for supervised classification, nonlinear SVMs stand out for providing robust and accurate decision boundaries, even in complex datasets {\cite{VapnikSVM1995}}. %
{Existing works that impose hard‐cardinality constraints to enforce strict feature selection mainly focus on the primal formulation of linear SVMs, solving the resulting mixed‐integer problem in the original feature space. When embedding %
feature selection
in nonlinear kernels, however,  most approaches typically resort to soft‐regularization techniques or approximation schemes that do not guarantee an exact feature count (see Section~\ref{sec:related} below). In contrast, we develop a formulation that embeds a hard cardinality constraint on feature usage directly into the dual of a kernel SVM, thereby retaining strict control over sparsity while exploiting the full representational power of nonlinear kernels. We then exploit the decomposable structure of the resulting optimization problem to devise efficient combinatorial algorithms for its approximate solution}.

\subsection{Related work}\label{sec:related}
{We are interested in formulations that embed feature selection in the training model of the SVM classifier.}
The majority of the literature addresses linear SVMs, with most of those works relying on regularization approaches (see~\cite{Bradley1998FeatureSV,FungFSNewton,weston2003use,hastie2004entire,nguyen2010optimal,yuan2010comparison,wang2006doubly,zou2005regularization}, for instance). However, we focus on feature selection models where a cardinality constraint explicitly limits the number of features to be used. These models can be preferable to regularization-based approaches in scenarios where the end user aims to impose a desired number of features, since regularization models often do not guarantee control over the number of selected features. Moreover, carefully tuning the penalization parameter to simultaneously achieve a specific  {number of features used} and good predictive performance can be particularly challenging \cite{bomze2025feature}. %

Among mixed-integer approaches for linear SVMs incorporating a cardinality {control on the features}, \cite{chan2007direct} replaced the explicit {cardinality}  constraint with a nonlinear condition involving $\ell_1$- and $\ell_2$-norms, admitting convex relaxations but yielding weaker solutions that may violate the {upper bound on} feature {selected} ~\cite{ghaddar2018high}. To address this, \cite{ghaddar2018high} proposed an alternative relaxation based on iterative parameter adjustment via bisection.
A bilevel optimization model was proposed by \cite{agor2019feature}, where binary variables at the upper level explicitly control feature selection, solved via a genetic algorithm. \cite{maldonado2014feature} and \cite{labbe2019mixed} formulated Mixed-Integer Linear Programming models for incorporating a cardinality constraint, with \cite{labbe2019mixed} also considered margin maximization through a $\ell_1$ norm regularization and proposed heuristic and exact approaches for their solution. In \cite{aytug2015feature}, a feature selection model based on binary masks in the $\ell_2$-regularized SVM is proposed, solved through Generalized Benders Decomposition, although convergence issues were observed for larger datasets. Lastly, in~\cite{bomze2025feature}, the authors address a similar feature selection problem for $\ell_2$-regularized SVMs under a cardinality constraint. They propose a mixed-integer formulation with complementarity constraints and develop heuristic and exact solution approaches based on scalable semidefinite programming relaxations to efficiently solve the resulting problems.

Handling highly nonlinear input data requires methods beyond linear SVMs. {However,} feature‐selection techniques for nonlinear SVMs, especially mixed‐integer formulations, have been less explored in the literature due to their substantially greater computational complexity. Note that we are interested in selecting the most informative features from the input space, and not from the projected spaces (as specified in \cite{weston2003use}). Most embedded methods for nonlinear SVMs are based on regularization approaches, where a penalization term is added to the objective function to seek a trade-off between the standard classification objective and the complexity of the classifier. Among the ones using continous variables only, in \cite{neumann2005combined}, four novel continuous regularization-based approaches are proposed, solved using difference-of-convex functions programming. In \cite{Maldonado2011KPSVM}, an approximation of the $\ell_0$-norm of the feature vector is added as a penalty term in the objective function of the dual SVM formulation, and a solution strategy is proposed that requires tuning six hyperparameters. Another continuous regularization approach was proposed by \cite{allen2013automatic}, introducing an $\ell_1$-norm penalization term; the resulting model is solved using linearization strategies and gradient-based techniques within an alternating optimization framework. More recently, \cite{Cordero2021MMFSSVM} formulated the embedded feature selection problem as a min–max optimization, reformulated it into a single-level nonlinear model using duality theory, and proposed a solution strategy to find a local optimum.

Among embedded methods for nonlinear SVMs involving binary variables, ~\cite{weston2000feature} presents a model that minimizes a radius-margin bound on the leave-one-out error of the hard-margin SVM. The binary variables are relaxed to solve a penalized version of the proposed optimization problem via gradient-based approaches.
In \cite{MangNLFSSVM2007} %
the authors propose a mixed-integer model and algorithm that generate a nonlinear kernel-based classifier with reduced input space features. However, they considered a different classifier, based on the less-known generalised SVM paradigm introduced in \cite{mangasarian1998generalized}.
They develop an algorithmic procedure to find a stationary point for their optimization problem. While these two papers model FS for nonlinear SVMs using binary variables, both approaches depart from the standard dual SVM formulation.

\gnote{We note that a classical way to obtain an SVM with strict cardinality control on the number of
features is \emph{recursive feature elimination} (RFE)~\cite{guyon2002gene}, a greedy
backward-elimination procedure originally devised for linear SVMs and later extended to the nonlinear
case. Starting from the full set of features, RFE repeatedly trains an SVM on a
recursively reduced subset of features, at each iteration discarding the least relevant one according to
a ranking criterion derived from the trained model, and it can be stopped as soon as the desired number
of features is reached. In this paper we show how this can be regarded as a heuristic procedure to solve
our cardinality-constrained nonlinear SVM problem.}

\subsection{Contribution}

Unlike mixed-integer approaches for linear SVMs, which rely on the primal formulation and typically associate binary variables with forcing hyperplane coefficients to zero for unselected features, we focus on nonlinear SVMs and propose a Mixed-Integer Nonlinear Programming (MINLP) model where binary variables directly \emph{mask} the contribution of each feature, as in \cite{MangNLFSSVM2007}. Feature selection is enforced exactly through a cardinality constraint on these binary variables, a setting that has not been previously explored in the literature.

For general nonlinear kernels, we first propose a \emph{local search} metaheuristic that directly operates on the binary variables. This approach exploits the decomposable structure of the problem to explore the search space effectively and provide high-quality solutions. Building on this idea, we then devise an alternating optimization framework that updates the continuous and binary variables in turn. When the binary variables are fixed, the resulting problem reduces to a standard SVM, which can be efficiently solved using dedicated algorithms. Conversely, when the continuous variables are fixed, the subproblem in the binary variables is highly nonlinear. For polynomial kernels, we show that this subproblem can be reformulated as the classical problem of maximizing a submodular function under a cardinality constraint.  This reformulation allows us to exploit scalable heuristics for submodular function optimization\lnote{, yielding our \emph{alternating convex-submodular optimization} scheme}. %

Submodular optimization is often used in feature selection problems, since one often aims to select features that maximize, for instance, submodular functions such as mutual information and conditional entropy drop ~\cite{bao2022submodular,krause2007near}, facility location and coverage~\cite{liu2013submodular,wei2013using}.  To the best of our knowledge, our paper is the first to employ submodular optimization in the context of interpretable nonlinear SVMs.

Computational experiments confirm the effectiveness of our approaches compared to standard solution methods for the resulting MINLP.

\bigskip

{The remainder of this paper is structured as follows: Section~\ref{sec: 2} presents the basics of SVMs and introduces our MINLP formulation. In Section~\ref{sec: 3}, we develop local search metaheuristics based on the decomposable structure of the MINLP. Section~\ref{sec: 4} shows that, for polynomial kernels, fixing the continuous variables yields a subproblem that can be reformulated as the maximization of a submodular function under cardinality constraints; based on this result, we propose a decomposition algorithm. Section~\ref{sec: 5} reports the numerical experiments, and Section~\ref{sec: 6} concludes the paper.}

\section{Embedded feature selection for nonlinear SVMs}
\label{sec: 2}

In this section we first recall the basics on nonlinear SVMs for classification and then describe the mixed-integer nonlinear formulation for embedded feature selection.

\subsection{Nonlinear Support Vector Machines}

Given a dataset $\{(x^i,y^i) \in \mathbb{R}^n\times\{-1,1\},\ i= 1, \dots, m\}$, we say that $x^i$ is the \textit{feature vector} of sample $i$ and $y^i$ its \textit{label}. Support Vector Machines (SVMs)~\cite{vapnik1999nature} are supervised machine learning models that seek a separating function with maximum margin.

In the  linear SVM formulation, the classifier is defined as the function $\fnote{\mathrm{Class}}:\mathbb{R}^n\to\{-1,1\}$ given by
\begin{equation}
    \fnote{\mathrm{Class}}(x) = \sgn({w^*}^T x + b^*),
\end{equation}
where the tuple $(w^*,b^*) \in \mathbb{R}^n \times \mathbb{R}$ is obtained as the optimal solution of the following convex quadratic optimization problem:
\begin{equation}
\label{prob: primal svm}
\begin{array}{rl}
    \min\limits_{w,b,\xi} \quad & \dfrac{1}{2} \|w\|_2^2 + C \sum_{i=1}^m \xi_i \\
    \text{s.t.} \quad
    & y^i(w^T x^i + b) \ge 1 - \xi_i \quad \forall i=1,\dots,m, \\
    & \xi_i \ge 0 \quad \forall i=1,\dots,m.
\end{array}
\end{equation}

When data is not linearly separable in the input space, the feature vectors are implicitly mapped into a higher-dimensional (possibly infinite-dimensional) space $\mathcal{H}$ through an implicit nonlinear transformation $\phi:\mathbb{R}^n\to\mathcal{H}$ and a linear SVM is applied in the transformed space. %
When nonlinear SVM are addressed, the dual formulation \cite{6789709}  of \eqref{prob: primal svm} is used:

\begin{equation}
\label{prob: dual svm}
\begin{array}{rl}
    \min\limits_{\alpha} \quad & \dfrac{1}{2} \alpha^T Q \alpha - \mathbf{1}^T \alpha \\
    \text{s.t.} \quad
    & y^T\alpha = 0, \\
    & 0 \leq \alpha \leq C,
\end{array}
\end{equation}
where $Q\in\mathbb{R}^{m\times m}$ is the matrix with entries
$$
q_{ih} = y^i y^h k(x^i, x^h).
$$
where the \textit{kernel function} $k(\cdot,\cdot)$ is defined as
$k(x^i,x^h) = \langle \phi(x^i), \phi(x^h) \rangle,$

Common examples of kernel functions include polynomial kernels, $k(x,z) = (\gamma x^Tz + c)^d$, or Gaussian ones, $k(x,z) = \exp(-\gamma \|x-z\|^2)$, where $ c \geq 0$, $\gamma > 0$  and $d \geq 1$.
Since kernel functions are defined to generate symmetric positive semidefinite matrices, matrix $Q$ is symmetric and positive semidefinite, ensuring the convexity of the dual problem~\cite{shawe2004kernel}.

\textcolor{revblue}{
Given an optimal solution $\alpha^*$ to problem~\eqref{prob: dual svm}, the classifier can be expressed as
\begin{equation*}
    \fnote{\mathrm{Class}}(x) = \sgn\left(\sum_{i=1}^m \alpha^*_i y^i k(x^i, x) + b^*\right),
\end{equation*}
where $b^*$ can be computed using standard reconstruction techniques (see~\cite{VapnikSVM1995}).}

\subsection{A mixed-integer nonlinear formulation}

Standard mixed-integer approaches, such as~\cite{labbe2019mixed, aytug2015feature, bomze2025feature}, based on the linear SVM formulation~\eqref{prob: primal svm}, model feature selection through a cardinality constraint and appropriate binary variables associated with the coefficients $w_j$. However, this way of modeling feature selection is not directly applicable when general nonlinear transformations of the input data are considered, because in nonlinear SVMs the mapping $\phi(\cdot)$ is not computed explicitly, and only scalar products evaluated via the kernel function are available. For this reason, we start from the dual formulation~\eqref{prob: dual svm} and, similarly to \cite{MangNLFSSVM2007, aytug2015feature}, we introduce binary variables $\beta \in \{0,1\}^n$ that mask  the contribution to the kernel function evaluation of original input features. The resulting problem is a mixed-integer nonlinear (and nonconvex) optimization problem.

\begin{equation}\label{prob:DFS-SVM}
\begin{array}{rl}
    \quad \min \limits_{\alpha, \beta}\quad & \displaystyle \frac{1}{2} \sum_{i=1}^{m} \sum_{h=1}^{m} y^i y^h \alpha_i \alpha_h k_{\beta}(x^i, x^h) - \sum_{i=1}^{m} \alpha_i \\
    \text{s.t.}  \quad & \displaystyle \sum_{i=1}^{m} y^i \alpha_i = 0 \\
    & \displaystyle \sum_{j=1}^{n} \beta_j = B \\
    & 0 \le \alpha_i \le C, \quad \forall i = 1, \dots, m \\
    & \beta_j \in \{0,1\}, \quad \forall j = 1, \dots, n \\
\end{array}
\end{equation}

where $$k_{\beta}(x^i,x^h) := k(\beta \odot x^i,\beta \odot x^h) = \langle \phi(\beta \odot x^i), \phi(\beta \odot x^h)\rangle$$ and $\odot$ denotes the element-wise (Hadamard) product between the vectors $\beta$ and $x^i$, and $B$ is a user-defined parameter. Considering that $\beta_j$ variables are binary, we have that:

\begin{itemize}
    \item for polynomial kernel functions, $ k_{\beta}(x^i, x^h) = \left( \displaystyle\gamma \sum_{j=1}^n {\beta_j} x^i_j  {\beta_j} x^h_j + c \right)^d  = \left( \displaystyle\gamma \sum_{j=1}^n {\beta_j} x^i_j x^h_j + c \right)^d;$
    \item for Gaussian kernel functions, $  k_{\beta}(x^i, x^h) =  \text{exp}\left ({-\gamma \ds \sum_{j=1}^n \left( {\beta_j} x^i_j -  {\beta_j}x^h_j \right)^2 } \right ) =  \text{exp}\left (-\gamma \ds \sum_{j=1}^n {\beta_j}\left(  x^i_j -  x^h_j \right)^2 \right ),$
\end{itemize}

{where $ c \geq 0$, $\gamma > 0$  and $d \geq 1$.} In this way, binary variables $\beta$ model the contribution of each feature in the kernel: if $\beta_j=1$, feature $j$ is included, otherwise it is excluded.

{
While solving problem~\eqref{prob:DFS-SVM} is challenging, it becomes more tractable when either the continuous or the binary variables are fixed. Indeed, problem~\eqref{prob:DFS-SVM} can be written as
$$
\begin{array}{rl}
    \displaystyle \min_{\alpha\in \R^m,\, \beta\in\{0,1\}^n} &  f(\alpha,\beta) =
    \displaystyle \frac{1}{2} \alpha^T Q(\beta)\alpha - e^T \alpha \\
    \text{s.t.}  \quad &
    (\alpha, \beta) \in \Delta \times \Omega,
\end{array}
$$
where the $m \times m$ symmetric matrix $Q(\beta)=\{y^i y^h k(\beta \odot x^i, \beta \odot x^h)\}_{i,h=1,\dots, m}$ is positive semidefinite for any $\beta$, and the two feasible regions are defined as
$\Delta=\{\alpha \in \R^{\textcolor{revblue}{m}}: \displaystyle\sum_{i=1}^m y^i \alpha_i = 0\textcolor{revblue}{,\ 0 \le \alpha_i \le C\ \,\forall i=1,\dots,m}\}$ and $\Omega=\{\beta \in \{0,1\}^n: \displaystyle\sum_{j=1}^n \beta_j = B\}$.

\lpnote{
\begin{remark}\label{rem:no-soft-reg}In the case of nonlinear SVM model \eqref{prob:DFS-SVM}, the
hard cardinality constraint $\sum_{j=1}^n\beta_j=B$ cannot be replaced by a soft $\ell_1$ regularization as usually done in linear SVM, see e.g. \cite{bomze2025feature}, without trivializing the model. In fact, if a penalty $\lambda\sum_{j=1}^n\beta_j$ (with $\lambda>0$) is added to $f(\alpha,\beta)$, the resulting problem admits a degenerate global optimum.
Indeed, by feasibility $\sum_{i\in I^+}\alpha_i=\sum_{i\in I^-}\alpha_i$ (with $I^{\pm}=\{i:y^i=\pm1\}$), and each side is at most $C\min\{|I^+|,|I^-|\}$, so we deduce $e^\top\alpha \le 2C\min\{|I^+|,|I^-|\}$.
Further, since $Q(\beta)\succeq0$ for every $\beta$, we have $\displaystyle f(\alpha,\beta)+\lambda\sum_{j=1}^n\beta_j\ge -e^\top\alpha \ge -2C\min\{|I^+|,|I^-|\}$, a value attained at $\beta^\star= 0$. Thus, a soft-regularized counterpart of model \eqref{prob:DFS-SVM} would simply select no features; it is precisely the hard cardinality constraint that enforces a meaningful, non-degenerate selection.
\end{remark}}

For a fixed binary vector $\beta = \bar \beta \in \{0,1\}^n$, the resulting problem, that we will refer to as the $\alpha$-subproblem, reduces to the dual SVM formulation \eqref{prob: dual svm} restricted to the features selected by $\bar\beta$, with kernel matrix $Q = \{q_{ih}\}_{i,h=1,\dots,m}$ defined as
$$
q_{ih} = y^i y^h k_{\bar \beta}(x^i, x^h).
$$

Conversely, for any fixed vector $\bar\alpha \in \R^m$, problem~\eqref{prob:DFS-SVM} becomes the following nonlinear binary problem, that we name the $\beta$-subproblem:
\begin{equation}\label{model:DFS-SVM-beta}
\begin{array}{rl}
    \min \limits_{\beta} \quad & \frac{1}{2}\displaystyle \sum_{i=1}^{m}\sum_{h=1}^{m} \bar{\alpha}_i \bar{\alpha}_h y^i y^h k_{\beta}(x^i, x^h) \\
    \text{s.t.} \quad & \displaystyle \sum_{j=1}^{n} \beta_j = B \\
    & \beta_j \in \{0,1\}, \quad \forall j = 1, \dots, n \\
\end{array}
\end{equation}

Our approaches leverage this decomposable structure by alternating optimization in the $\alpha$ and $\beta$ variables. When $\beta = \bar\beta$ is fixed, the convex problem in $\alpha$ can be efficiently solved, \gnote{for instance by the Sequential Minimal Optimization (SMO) algorithm}, as implemented in \texttt{scikit-learn}~\cite{scikit-learn}. In contrast, optimizing over $\beta$ leads to a non-convex binary problem, particularly difficult in the presence of nonlinear kernels. Moreover, a naive alternating approach iteratively solving the $\alpha$-subproblem and the $\beta$-subproblem does not guarantee convergence to the global optimum of~\eqref{prob:DFS-SVM} and often yields poor-quality heuristic solutions in practice.

To address this, we propose:
\begin{itemize}
    \item A local search framework for general kernels, which updates $\beta$ by exploring small neighborhoods while repeatedly solving the $\alpha$-subproblem. An improved variant samples from larger neighborhoods to escape local minima;

    \item A submodular reformulation of the $\beta$-subproblem~\eqref{model:DFS-SVM-beta} for polynomial kernels, enabling the integration of scalable greedy algorithms in a more effective alternating optimization scheme.

\end{itemize}}

\section{Local Search algorithms for general kernels \label{sec: 3}}

We first introduce a local search algorithm that is able to find a local solution to problem \eqref{prob:DFS-SVM} for general kernels; this framework is then integrated into the subsequent algorithms presented in the paper. We then propose an improved version of the local search algorithm with enhanced exploration capabilities of the feasible region. It is worth noting that both algorithms do not tackle the full MINLP formulation directly. Instead, they sidestep the joint optimization by iteratively fixing the binary variables $\beta$ and solving the resulting continuous SVM dual subproblem. \gnote{In this way they never face the nonconvexity of~\eqref{prob:DFS-SVM} as a whole, and the entire search effort is devoted to the binary space.}

\subsection{\gnote{A Local Search algorithm}}

Let us introduce the following neighborhood, denoted by $\mathcal{N}^1(\beta^k)$, for a given feasible solution $\beta^k$ to \eqref{model:DFS-SVM-beta}:
$$
\mathcal{N}^1(\beta^k) := \left\{ \beta \in \{0,1\}^n : \sum_{j \in [n]} \beta_j = B,\; \gnote{s(\beta,\beta^k) = 1} \right\},
$$
\gnote{where
$$
s(\beta,\beta^k) := \bigl| \{\, j \in [n] \;:\; \beta^k_j = 1,\; \beta_j = 0 \,\} \bigr|
$$
is the number of features selected by $\beta^k$ that $\beta$ discards.} The set $\mathcal{N}^1(\beta^k)$ thus consists of binary vectors $\beta$ that differ from $\beta^k$ by flipping one entry from 1 to 0 and one from 0 to 1\gnote{, that is, by a single swap of a selected feature for an unselected one}, thereby preserving the cardinality constraint $\sum_{j \in [n]} \beta_j = B$. The size of $\mathcal{N}^1(\beta^k)$ is $|\mathcal{N}^1(\beta^k)| = B(n-B) = nB - B^2$.

Our first local search algorithm starts from a feasible solution $\beta^{\text{in}}$, which can be generated either randomly or via a heuristic, and sets it as the current solution $\beta^k$. At each iteration, the algorithm explores the neighborhood $\mathcal{N}^1(\beta^k)$ \lpnote{of cardinality $B(n-B)$}. For each $\widehat \beta \in \mathcal{N}^1(\beta^k)$, it solves the corresponding dual SVM problem, obtaining the objective value of the dual formulation (denoted as \lpnote{$\widehat{\text{UB}}$} ) and the associated dual variables $\widehat \alpha$. The algorithm keeps track of the best objective value found so far, denoted UB$^*$, and updates the best solution whenever an improvement is achieved. Once the entire neighborhood has been explored, the algorithm moves to the best solution found within it and repeats the process until no further improvement is observed. \textcolor{revblue}{We note that the $B(n-B)$ dual SVM problems associated with the neighbours in $\mathcal{N}^1(\beta^k)$ are mutually independent because each simply fixes $\beta=\widehat\beta$ and solves the corresponding convex SVM. Thus, the neighborhood exploration can be carried out in parallel. We detail how we exploit this parallelism in Appendix~\ref{app:param}.} We denote this algorithm by \lpnote{\emph{Local Search}} \gnote{(}\texttt{LS}\gnote{)}; see Algorithm~\ref{alg:LS}.

\begin{algorithm}
\caption{\texttt{Local Search} (\texttt{LS})}
\label{alg:LS}
\begin{algorithmic}[1]
\State \textbf{Input:} Initial solution $\beta^{\text{in}} \in \{0,1\}^n$
\State \textbf{Output:} Solution $(\alpha^*, \beta^*)$ to \eqref{prob:DFS-SVM} and value UB$^*$
\State Fix $\beta = \beta^{\text{in}}$, solve \eqref{prob: dual svm} to obtain $(\alpha^{\text{in}}, \text{UB}^{\text{in}})$
\State Initialize $k \gets 1$, $\beta^{k-1} \gets \mathbf{0}$, $\beta^k \gets \beta^{\text{in}}$, $\alpha^k \gets \alpha^{\text{in}}$,  UB$^* \gets \text{UB}^{\text{in}}$
\While{$\beta^k \neq \beta^{k-1}$}
    \State \textcolor{revblue}{Set $(\alpha^{k+1}, \beta^{k+1}) \gets (\alpha^{k}, \beta^{k})$}
    \For{each $\widehat{\beta} \in \mathcal{N}^1(\beta^k)$}
        \State Fix $\beta = \widehat{\beta}$ and solve \eqref{prob: dual svm} to obtain $(\widehat{\alpha}, \widehat{\text{UB}})$
        \If{$\widehat{\text{UB}} \textcolor{revblue}{<} \text{UB}^*$}
            \State Set $(\alpha^{k+1}, \beta^{k+1}, \text{UB}^*) \gets (\widehat{\alpha}, \widehat{\beta}, \widehat{\text{UB}})$
        \EndIf
    \EndFor
    \State Update $k \gets k+1$
\EndWhile
\State Set $(\alpha^*, \beta^*) \gets (\alpha^k, \beta^k)$
\State \textbf{return} $(\alpha^*, \beta^*, \text{UB}^*)$
\end{algorithmic}
\end{algorithm}

In practice, this algorithm is quite fast on standard datasets, such as those used in our computational experiments. Nevertheless, it can be improved to find better solutions.

\subsection{An improved Local Search metaheuristic}

We now introduce a broader neighborhood, denoted by $\mathcal{N}^2(\beta^k)$. Specifically, we define:
$$
\mathcal{N}^2(\beta^k) := \left\{ \beta \in \{0,1\}^n : \sum_{j \in [n]} \beta_j = B,\; \gnote{s(\beta,\beta^k) \geq 2} \right\},
$$
\gnote{that is, $\mathcal{N}^2(\beta^k)$ contains the solutions obtained from $\beta^k$ by swapping at
least two selected features at once, always preserving the
cardinality constraint $\sum_{j \in [n]} \beta_j = B$. Single swaps are thus excluded, so
that $\mathcal{N}^1(\beta^k)$ and $\mathcal{N}^2(\beta^k)$ are disjoint.}
$\mathcal{N}^2(\beta^k)$ is significantly larger than $\mathcal{N}^1(\beta^k)$, but it can be explored by randomly sampling binary vectors $\beta$ from within $\mathcal{N}^2(\beta^k)$\gnote{, drawing for each sample the number of swapped couples uniformly among the admissible values}.

The improved version of the local search algorithm proceeds as follows. After running \texttt{LS} starting from an initial solution $\beta^k_1$, the algorithm obtains a new solution $\beta^k_2$. To explore other promising regions of the feasible space, the algorithm samples the neighborhood $\mathcal{N}^2(\beta_2^k)$, generating \texttt{$\#$Samples} candidate vectors $\widehat{\beta}$. For each of them, the corresponding dual SVM problem is solved, and the best solution is updated if an improvement is found. The algorithm then moves to the best $\widehat{\beta}$ identified within this sample and restarts the  \texttt{LS} procedure.

\lnote{A further enhancement of \texttt{LS$^*$} reduces the cost of the inner local search over $\mathcal{N}^1$: rather than re-fitting the dual variables for every neighbour, \texttt{LS$^*$} \emph{screens} all candidates in $\mathcal{N}^1(\beta^k)$ at the current dual solution $\alpha^k$, evaluating the objective of \eqref{prob: dual svm} at $(\alpha^k,\widehat\beta)$ without re-solving the SVM, and re-solves the SVM (the exact re-fit) only for the top fraction, the \emph{re-fit fraction} $\%$re-fit, retaining the best neighbour by the true objective. Re-fitting all neighbours ($\%\mathrm{re\text{-}fit}=100\%$) recovers the plain \texttt{LS}.}

To make sure that the same neighborhood $\mathcal{N}^1(\beta^k)$ is never explored more than once, we introduce a tabu list to store all solutions $\beta^k$ whose $\mathcal{N}^1(\beta^k)$ neighborhood has already been visited. This list is used to avoid redundant exploration and can be integrated into \texttt{LS}. As a result, the algorithm may sometimes accept a new $\beta$ that leads to a worse objective function value; however, this new $\beta$ serves as a fresh starting point for \texttt{LS}, possibly guiding the search to a better region of the feasible set. \textcolor{revblue}{We remark that, although the next starting point $\Tilde{\beta}$ is selected as the sampled candidate with the best objective value, this choice does not conflict with the diversification purpose of this phase. Indeed, the candidates are drawn from the enlarged neighborhood $\mathcal{N}^2(\beta_2^k)$, which contains \gnote{the} multi-swap moves that the single-swap neighborhood $\mathcal{N}^1$ explored by \texttt{LS} cannot reach; moreover, the tabu list forbids any already-visited $\mathcal{N}^1$ neighborhood. Hence $\Tilde{\beta}$ always seeds \texttt{LS} in a region genuinely different from the one just explored, and the objective value is used only as an inexpensive proxy to pick the most promising such region: the escape from the current local optimum is guaranteed by the sampling of $\mathcal{N}^2$ and by the tabu mechanism, not by the selection rule itself.}

The algorithm terminates when a stopping criterion is met. Possible criteria include a maximum number of iterations, a CPU time limit, or a maximum number of consecutive iterations without improvement. In our implementation, we adopt the last option, and we denote by \gnote{\texttt{OptWindow} the corresponding parameter, that is, the number of consecutive non-improving iterations tolerated before stopping}. We refer to this improved algorithm as \gnote{\emph{Improved Local Search} (}\texttt{LS$^*$}\gnote{)}; see Algorithm~\ref{alg:LS_star}.

\begin{algorithm}
\caption{\texttt{Improved Local Search} (\texttt{LS$^*$})}
\label{alg:LS_star}
\begin{algorithmic}[1]
\State \textbf{Input:} Initial solution $\beta^{\text{in}} \in \{0,1\}^n$, number of samples \texttt{\#Samples}\lnote{, re-fit fraction $\%$re-fit}
\State \textbf{Output:} Solution $(\alpha^*, \beta^*)$ to \eqref{prob:DFS-SVM} and value UB$^*$
\State Initialize $k \gets 1$, $\beta_1^{k} \gets \beta^{\text{in}}$, \texttt{TabuList} $\gets \emptyset$, UB$^* \gets \infty$

\While{termination criteria not met}
    \State Perform \texttt{LS} \lnote{(with re-fit fraction $\%$re-fit)} starting from $\beta_1^k$ and populate \texttt{TabuList} with solutions whose $\mathcal{N}^1$ has been explored
    \State Set $(\alpha^k, \beta_2^k, \text{UB})$ as the best solution found by \texttt{LS}
    \If{$\text{UB} \leq \text{UB}^*$}
        \State Update $(\alpha^*, \beta^*, \text{UB}^*) \gets (\alpha^k, \beta_2^k, \text{UB})$
    \EndIf
    \State Sample $\mathcal{N}^2(\beta_2^k)$ to obtain \texttt{\#Samples} candidate solutions
    \For{each sampled $\widehat{\beta} \in \mathcal{N}^2(\beta_2^k)$ not in \texttt{TabuList}}
        \State Fix $\beta = \widehat{\beta}$ and solve \eqref{prob: dual svm} to obtain $(\widehat{\alpha}, \widehat{\text{UB}})$
        \If{$\widehat{\text{UB}} \leq \text{UB}^*$}
            \State Update $(\alpha^*, \beta^*, \text{UB}^*) \gets (\widehat{\alpha}, \widehat{\beta}, \widehat{\text{UB}})$
        \EndIf
    \EndFor
    \State Set $\Tilde{\beta}$ as the best sampled $\widehat{\beta}$ in $\mathcal{N}^2(\beta_2^k)$
    \State Update $\beta_1^{k+1} \gets \Tilde{\beta}$
    \State Update $k \gets k+1$
\EndWhile
\State \textbf{return} $(\alpha^*, \beta^*, \text{UB}^*)$
\end{algorithmic}
\end{algorithm}

\section{\gnote{An alternating optimization} algorithm for polynomial kernels}
\label{sec: 4}

When polynomial kernels are used, a direct approach to solving the $\beta$-subproblems~\eqref{model:DFS-SVM-beta} relies on linearized reformulations~\cite{mythesis2024}. %
However, as both the number of features $n$ and the polynomial degree $d$ increase, their size grows rapidly, making them computationally burdensome.

In this section, we propose an alternative approach that reformulates subproblem~\eqref{model:DFS-SVM-beta} as the maximization of a submodular set function. Unlike linearized formulations, the submodular reformulation decouples the complexity of the problem from the polynomial degree $d$, leading to a more scalable and efficient method. We then show how this reformulation can be integrated into a new alternating optimization scheme.

\subsection{Submodular functions optimization}
\label{sec: sub fun opt}

Let $E$ be a finite set of elements. A set function $f: 2^E \rightarrow \mathbb{R}$ is called \textit{submodular} if for all subsets $S, T \subseteq E$, the following condition holds:
\begin{equation}
f(S) + f(T) \geq f(S\cup T)+ f(S \cap T).
\label{def: submodularity}
\end{equation}

Equivalently, a set function $f$ is submodular if and only if it satisfies the following property:
\begin{equation}
f(T \cup \{e\}) - f(T) \geq f(S \cup \{e\}) - f(S) \quad \forall \,  T \subseteq S \subseteq E,\; e \in E \setminus S.
\label{def: dr submodularity}
\end{equation}

For \emph{monotone non-decreasing} functions, the latter can be interpreted as a \emph{diminishing returns} property: adding an element to a smaller set contributes more to the function's increase than adding the same element to a larger set. On the contrary $f$ is said to be \textit{supermodular} if the sign of \eqref{def: submodularity} and \eqref{def: dr submodularity} is reversed; {equivalently, $f$ is supermodular if $-f$ is submodular.} %

One of the most studied problem involving {monotone non-decreasing} submodular functions is their cardinality constrained maximization problem:
\begin{align}
\max_S \; & f(S) \label{model: max submod} \\
\text{s.t.}\;  & S \subseteq E \nonumber\\
& |S| \leq k \nonumber
\end{align}
for a given $0 \leq k \leq |E|$ and $f$ submodular {and monotone non-decreasing}. This optimization problem is in general $\mathcal{NP}$-hard \cite{nemhauser1978best,}. Nevertheless, it is well-known that a simple greedy algorithm achieves an approximation ratio of $(1 - \frac{1}{e})$ {in theory~\cite{nemhauser1978analysis}  (and often better in practice \cite{krause2014submodular})} under the additional hypothesis that $f$ is nonnegative.
The simple greedy can be improved in practice by exploiting definition \eqref{def: dr submodularity}, i.e., the fact that as we add elements to the set, the marginal benefits of the remaining elements to consider can only decrease. {This observation lead to the \texttt{Lazy Greedy} algorithm, introduced} by Minoux \cite{minoux2005accelerated}, {which in practice} reduces the number of evaluations of the submodular function, leading to faster convergence, especially in large-scale problems.

\textcolor{revblue}{Concretely, \texttt{Lazy Greedy} maintains, for every element $e\in E$, an upper bound on its marginal gain $\Delta_e(S)=f(S\cup\{e\})-f(S)$ and keeps these bounds in a max-priority queue. At each iteration it pops the top element and recomputes its marginal gain with respect to the current set $S$: by submodularity (the diminishing-returns property~\eqref{def: dr submodularity}) this gain can only have decreased, so if the recomputed value still dominates the next bound in the queue, that element is added to $S$ without evaluating any other; otherwise its updated bound is re-inserted and the next element is examined. This lazy re-evaluation typically avoids recomputing the vast majority of marginal gains, while returning exactly the same set as the standard greedy.}

\subsection{A reformulation based on submodularity}

Here we show that it is possible to formulate problem \eqref{model:DFS-SVM-beta} as \eqref{model: max submod}, thus as a maximization of a {nonnegative, monotonically non-decreasing} submodular set function problem, {under a cardinality constraint}. Let $N= 2^{[n]} $ be the power set of $[n]$, i.e., the set that contains all possible subsets of the set $[n]$.  For a subset of features $S \in N$ let us define  $x^i_S$ as the feature vector of sample $i$ indexed only over the features in $S$. Given a kernel function $k(\cdot,\cdot)$, we define the kernel matrix as the matrix $K_S$ with components $\left[K_S \right]_{ih} = k(x^i_S, x^h_S)$, $i,h\in [m]$. Moreover let us define the set function $F(S)$ as $F(S)=0$ for $S=\emptyset$ and
$$F(S) = \sum_{i=1}^m\sum_{h=1}^m \bar\alpha_i\bar \alpha_h y^iy^h k(x^i_S, x^h_S),$$
otherwise.

Note that, since $ k(\cdot, \cdot)$ is a kernel function, for any set $S \in N$, the matrix $Q_S \in \R^{m \times m}$ with components $q_{ih} = y^i y^h k(x_S^i, x_S^h)$, is positive semidefinite. Thus for any fixed vector $\bar \alpha$, $F(S)$ is always non negative:
\begin{equation}\label{eq: nonegativity of F}
F(S) = \bar \alpha^T Q_S \bar \alpha \geq 0, \quad \forall S\in N.
\end{equation}

We now restrict to the case in which polynomial kernels are adopted. Thus, we consider set functions of the type:
\begin{align*}
F(S) = \sum_{i=1}^m\sum_{h=1}^m \bar\alpha_i\bar \alpha_h y^iy^h \left (c + \sum_{j \in S} \gamma x^i_j x^h_j \right )^d,  & \quad d \geq 2, S \in N\setminus \textcolor{revblue}{\{\emptyset\}}.
\end{align*}

Given this notation, we can first rewrite problem \eqref{model:DFS-SVM-beta}  as the following minimization problem:
\begin{equation}\label{model: min supermod}
\begin{array}{rl}
    \min\limits_{S} & F(S) \\
    \text{s.t.} & S \in  N \\
    & |S| = B.
\end{array}
\end{equation}

The following proposition, whose proof can be found in \gnote{Appendix~\ref{app:proofs}}, states that the set function $F$ is a monotone non-decreasing supermodular set function when polynomial kernel functions are adopted.

\begin{proposition}
    If $k(\cdot, \cdot)$ is a polynomial kernel function, the set function F(S) is %
    supermodular and monotone non-decreasing.
    \label{prop: submodularity and monotonicity of F}
\end{proposition}

Starting from this result, we can then reduce problem \eqref{model: min supermod} to a submodular set function maximization problem under cardinality constraint. Recall that the submodular objective function needs to be monotone nondecreasing and nonnegative. Consider the following problem:

\begin{equation}\label{model: reduction to max submod}
\begin{array}{rl}
\max \limits_S \; & H(S)  \\
\text{s.t.}\;  &  S \in N  \\
& |S| \leq n - B
\end{array}
\end{equation}

where $$H(S) = \lambda -  \sum_{i=1}^m\sum_{h=1}^m \bar\alpha_i\bar \alpha_h y^iy^h  \left (c + \sum_{j \in [n] \setminus  S} \gamma x^i_j x^h_j \right )^d   $$ and $\lambda =  \sum_{i=1}^m\sum_{h=1}^m \bar\alpha_i\bar \alpha_h y^iy^h \left (c + \sum_{j \in [n]} \gamma x^i_j x^h_j \right )^d$ .

We can now state our main result in the next proposition whose proof can also be found in \gnote{Appendix~\ref{app:proofs}}.

\begin{proposition}
    \label{prop: H sub}
    Function $H$ is a monotone nondecreasing and nonnegative submodular function. Problems \eqref{model: reduction to max submod} and \eqref{model: min supermod} are equivalent in the sense that, given an optimal solution $S''$ of \eqref{model: reduction to max submod} with value $H(S'')$, we can obtain an optimal solution $S'$ of \eqref{model: min supermod} as $S'= [n]\setminus S''$ with value $F(S')=\lambda - H(S'')$.
\end{proposition}

\subsection{An algorithm based on the submodular reformulation}

Depending on how the problem in the $\beta$ variables \eqref{model:DFS-SVM-beta} is solved\textcolor{revblue}{,} we can have different versions of the same decomposition framework. In this paper we show the decomposition algorithm that yielded the best performances in our computational experiments, which is the one where the problem in the $\beta$ variables is solved through its submodular reformulation with the \texttt{Lazy Greedy} algorithm by \cite{minoux2005accelerated}.  We invite the interested reader to consult \cite{mythesis2024}, where we show all the details and the experiments on the other decomposition algorithms solving the problem in the $\beta$ variables by either linearization procedures or by a MIP model for submodular function optimization.

\fnote{We first present the basic scheme, the \emph{Alternating Convex-Submodular Optimization} (\AObase). Starting from an initial feature subset $\beta^{\text{in}}$, it alternates between the two subproblems: for a fixed $\beta$ it solves the convex $\alpha$-subproblem~\eqref{prob: dual svm} with a standard SVM solver, and for the fixed $\alpha$ it solves the $\beta$-subproblem~\eqref{model:DFS-SVM-beta} through its submodular reformulation using the \texttt{Lazy Greedy} algorithm of \cite{minoux2005accelerated}; the procedure terminates at the first iteration that does not improve the objective value. The scheme is summarized in Algorithm~\ref{alg:ACSO}.}

\begin{algorithm}
\caption{\fnote{Alternating Convex-Submodular Optimization} (\AObase)}
\label{alg:ACSO}
\begin{algorithmic}[1]
\State \textbf{Input:} Initial solution $\beta^{\text{in}} \in \{0,1\}^n$
\State \textbf{Output:} Solution $(\alpha^*, \beta^*)$ to \eqref{prob:DFS-SVM} and value $\text{UB}^*$
\State Fix $\beta = \beta^{\text{in}}$, solve \eqref{prob: dual svm} to obtain $(\alpha^*, \text{UB}^*)$; set $\beta^* \gets \beta^{\text{in}}$
\State \texttt{improved} $\gets$ \textbf{true}
\While{\texttt{improved}}
    \State Fix $\alpha = \alpha^*$, solve the submodular maximization~\eqref{model: reduction to max submod} to obtain $\Tilde{\beta}$
    \State Fix $\beta = \Tilde{\beta}$ and solve \eqref{prob: dual svm} to obtain $(\Tilde{\alpha}, \widetilde{\text{UB}})$
    \If{$\widetilde{\text{UB}} < \text{UB}^*$}
        \State Update $(\alpha^*, \beta^*, \text{UB}^*) \gets (\Tilde{\alpha}, \Tilde{\beta}, \widetilde{\text{UB}})$
    \Else
        \State \texttt{improved} $\gets$ \textbf{false}
    \EndIf
\EndWhile
\State \textbf{return} $(\alpha^*, \beta^*, \text{UB}^*)$
\end{algorithmic}
\end{algorithm}

\fnote{Although simple, \AObase{} often gets trapped in poor-quality local optima. To overcome this limitation, we introduce an enhanced version, the \emph{Improved Alternating Convex-Submodular Optimization} (\AOsub), obtained from \AObase{} by adding two mechanisms.} First, a \emph{multisolution} strategy is adopted: for each fixed $\alpha$, multiple candidate $\beta$ solutions are generated and stored in a pool (\texttt{\#SolPool}); each candidate is then used to recompute $\alpha$, and the best pair is retained. \textcolor{revblue}{As in the local search, these \texttt{\#SolPool} re-computations of $\alpha$ are mutually independent (each fixes a candidate $\beta$ and solves the corresponding convex SVM) and can therefore be carried out in parallel.} Second, to prevent revisiting previously explored $\beta$ solutions, a tabu list is maintained and updated at each iteration. These mechanisms allow the algorithm to explore a broader region of the search space and improve the quality of the final solution compared to naive alternation.

\begin{algorithm}
\caption{\fnote{Improved Alternating Convex-Submodular Optimization} (\AOsub)}
\label{alg:DecSubLight$^*$}
\begin{algorithmic}[1]
\State \textbf{Input:} Initial solution $\beta^{\text{in}} \in \{0,1\}^n$, integer \texttt{\#SolPool}
\State \textbf{Output:} Solution $(\alpha^*, \beta^*)$ to \eqref{prob:DFS-SVM} and value UB$^*$
\State Fix $\beta = \beta^{\text{in}}$, solve \eqref{prob: dual svm} to obtain $(\alpha^{\text{in}}, \text{UB}^{\text{in}})$
\State Initialize $(\alpha^k, \beta^k) \gets (\alpha^{\text{in}}, \beta^{\text{in}})$, \texttt{TabuList} $\gets \emptyset$, \texttt{ExploredSolsList} $\gets \emptyset$, $\text{UB} \gets \text{UB}^{\text{in}}$

\While{termination criteria not met}
    \State Fix $\alpha = \alpha^k$, solve the submodular maximization~\eqref{model: reduction to max submod} to obtain solution $\Tilde{\beta}$
    \State Apply \texttt{LS} starting from $\Tilde{\beta}$ without recomputing $\alpha$ values
    \State Retrieve the best \texttt{\#SolPool} solutions $\widehat{\beta}$ found
    \For{each $\widehat{\beta}$ not in \texttt{TabuList}}
        \State Fix $\beta = \widehat{\beta}$ and solve \eqref{prob: dual svm} to obtain $(\widehat{\alpha}, \widehat{\text{UB}})$
        \State Add $\widehat{\beta}$ to \texttt{TabuList}
        \State Add $(\widehat{\alpha}, \widehat{\beta}, \widehat{\text{UB}})$ to \texttt{ExploredSolsList}
        \If{$\widehat{\text{UB}} \leq \text{UB}^*$}
            \State Update  $(\alpha^*, \beta^*, \text{UB}^*) \gets (\widehat{\alpha}, \widehat{\beta}, \widehat{\text{UB}})$
        \EndIf
    \EndFor
    \State Set $(\alpha^{k+1}, \beta^{k+1}) \gets (\textcolor{revblue}{\bar{\alpha}}, \textcolor{revblue}{\bar{\beta}})$ where $\textcolor{revblue}{(\bar{\alpha}, \bar{\beta}, \overline{\text{UB}})}$ is the best solution in \texttt{ExploredSolsList} \textcolor{revblue}{(i.e., the triple with the smallest $\overline{\text{UB}}$)}
    \State Remove $\textcolor{revblue}{(\bar{\alpha}, \bar{\beta}, \overline{\text{UB}})}$ from \texttt{ExploredSolsList}
    \State Update $k \gets k+1$
\EndWhile
\State \textbf{return} $(\alpha^*, \beta^*, \text{UB}^*)$
\end{algorithmic}
\end{algorithm}

Algorithm~\ref{alg:DecSubLight$^*$} illustrates the decomposition framework. The algorithm starts from an initial solution $\beta^{\text{in}}$ and, at the beginning of each iteration, fixes the continuous variables $\alpha$ to the current best solution $\alpha^k$. A candidate solution $\Tilde{\beta}$ is then computed using the \texttt{Lazy Greedy} algorithm. Starting from $\Tilde{\beta}$, the \texttt{LS} procedure is applied while keeping $\alpha$ fixed to $\alpha^k$, i.e., without recomputing $\alpha$ for each $\beta$ explored during the search. The best \texttt{\#SolPool} solutions $\widehat{\beta}$ found through this process are collected in the list \texttt{ExploredSolsList}, which stores candidate solutions to be evaluated in subsequent iterations.

For each $\widehat{\beta}$ in \texttt{ExploredSolsList} that has not been previously evaluated (tracked via a tabu list), the algorithm solves the dual SVM problem with $\beta = \widehat{\beta}$, thus recomputing the corresponding $\widehat{\alpha}$ and objective value ($\widehat{\text{UB}}$). If an improvement in UB$^*$ is observed, the best solution is updated accordingly. The algorithm then proceeds by selecting the best candidate from \texttt{ExploredSolsList} as the new current solution $(\alpha^{k+1}, \beta^{k+1})$ and repeats this process. The search terminates when a stopping criterion is met, such as a maximum number of iterations, time limit, or lack of improvement for a fixed number of consecutive iterations (the latter being used in our implementation, with the parameter \gnote{\texttt{OptWindow}}).

Note that in Algorithm \ref{alg:DecSubLight$^*$}, we solve the subproblems in the $\beta$ variables heuristically using \texttt{Lazy Greedy}. Another version, which can be found in \cite{mythesis2024}, solves each of these subproblems to optimality through the MIP formulation and approach proposed in \cite{salvagnin2019some}. However, this latter approach resulted in much longer computational times, and the final solutions found by the two decomposition approaches were always the same. This is likely due to two reasons: first, in many of our experiments, the \texttt{Lazy Greedy} solution was the same as the one found by the global approach; and second, because we are using the multisolution approach, when recomputing the best $\alpha$ for each $\beta$ in the pool, the algorithm often moves to suboptimal $\beta$ solutions, discarding the optimal one.

\section{Computational results}
\label{sec: 5}
\begin{table}[ht]
\centering
\renewcommand\arraystretch{1.2}
\begin{tabular}{lccc}
\toprule
Dataset & $m$ & $n$ & Class (\%) \\
\midrule
Wholesale              & 440  & 7   & 68/32 \\
Diabetes               & 768  & 8   & 65/35 \\
\gnote{Breast Cancer Wisconsin (BCW)} & 699  & 9   & {66/34} \\
Cleveland              & 297  & 13  & 54/46 \\
Parkinsons             & 195  & 22  & 25/75 \\
German                 & 1000 & 24  & 70/30 \\
\gnote{Breast Cancer Diagnostic (BCD)} & 569  & 30  & 63/37 \\
Ionosphere             & 351  & 33  & 36/64 \\
Sonar                  & 208  & 60  & 53/47 \\
\bottomrule
\end{tabular}
\caption{\textcolor{revblue}{Datasets used in the computational study: \lpnote{$n$ number of
features; $m$ number of samples; Class (\%) number of samples in the two classes.}}}
\label{table:small_datasets_1}

\hfill

\end{table}

In this section, we analyze our proposed algorithms for solving problem~\eqref{prob:DFS-SVM}.  The experiments were conducted on the datasets listed in Table~\ref{table:small_datasets_1}\textcolor{revblue}{, all of which are publicly available from the UCI Machine Learning Repository\footnote{\texttt{https://archive.ics.uci.edu/}}}. All computations were performed on an Apple M1 CPU and 16 GB of RAM computer \lpnote{except for the enumeration approach reported in Section \ref{sec:certified}}. All models and algorithms were implemented in \texttt{Python}.
\gnote{The code implementing the proposed algorithms, together with the datasets used in this study, is
publicly available at \texttt{https://github.com/f-donofrio/fs-nl-svm}.}

{\color{revblue} %

Throughout this section the SVM hyper-parameters are held fixed,
\lpnote{\[
  \boxed{\;C=10 \qquad \gamma=0.1\;}
\]}for every dataset, kernel and budget.
\lpnote{T}he purpose %
is to compare the
\emph{optimization} %
\lpnote{performance} of the different methods on a common set of instances, so the instances
themselves must not change from one method to another. The effect of $C$ and $\gamma$ on the choice of
the parameters of our heuristics is examined in Appendix~\ref{app:param}.
\gnote{Throughout this section we report the results obtained by the improved versions \texttt{LS$^*$} and \AOsub; a comparison with their basic counterparts \texttt{LS} and \AObase, over which they consistently improve, is reported in Appendix~\ref{app:compare}.}

The algorithm-specific parameters of our metaheuristics are set, throughout this section, to
\[
  \boxed{\;\texttt{LS$^*$}:\ 10\%\text{ re-fit},\ \#\mathrm{Samples}=200,\ \texttt{OptWindow}=10\;}
  \quad
  \boxed{\;\AOsub:\ \#\mathrm{SolPool}=100,\ \texttt{OptWindow}=10\;}
\]
selected by the parameter study reported in Appendix~\ref{app:param}.
These values were selected on a broad pool of mixed instances
 as the configuration giving the best quality--cost trade-off across the whole grid. This study \lpnote{also allows us to compare}  the heuristics. Details are reported in Appendix \ref{app:compare}.
 The heuristics reach the same objective in most cases and, where they differ, \AOsub attains the better solution quality while \texttt{LS$^*$} solves fewer SVM subproblems, both being fast in absolute terms. The two heuristics differ in scope as well, since \AOsub is restricted to polynomial kernels by its submodular reformulation, so on the Gaussian kernel \texttt{LS$^*$} is the only applicable option.

 Since both \texttt{LS$^*$} and \AOsub need an initial $\beta^{\text{in}}$ to start with, we consider two possible cases: \begin{itemize}
    \item deterministic: the same starting point is provided to all the algorithms, which are run once;
    \item random: several runs from different starting points are considered, and mean values of the results are reported.
\end{itemize}

For the deterministic case, we derive a simple heuristic to obtain the starting point. Starting from the solution of the SVM dual \eqref{prob: dual svm} with all features, we keep the corresponding $\alpha$ values fixed and build the initial subset $\beta^{\text{in}}$ greedily, removing one feature at a time, always the one \gnote{whose removal yields the smallest value of the dual objective}, until $B$ features remain
(for polynomial kernels we use the faster  \texttt{Lazy Greedy} algorithm by \cite{minoux2005accelerated}, see Section \ref{sec: sub fun opt}). We name this initialization \emph{$\alpha$-start.} }

\subsection{\textcolor{revblue}{Comparison with other solution approaches}}\label{sec:comparison}

\lpnote{Throughout this section we compare the proposed heuristics against the general-purpose nonlinear solvers
\texttt{Knitro} (version 14.2.0) and \texttt{Ipopt} (version 3.14), both applied to
formulation~\eqref{prob:DFS-SVM} modeled in \texttt{Pyomo}, and against two
 heuristics based on the \emph{recursive feature elimination} (RFE)~\cite{guyon2002gene}.
}

{\color{revblue}

The two solvers play different roles. \texttt{Knitro} is a commercial solver that handles general MINLPs,
hence it tackles formulation~\eqref{prob:DFS-SVM} directly, with the binary variables $\beta$. Being a
\emph{local} (not global) MINLP solver, it returns feasible solutions with no optimality certificate.
\texttt{Ipopt} is an open-source interior-point solver for continuous nonlinear programs and cannot handle
integrality at all: we therefore solve the continuous relaxation of~\eqref{prob:DFS-SVM}, with
$\beta\in[0,1]^n$, and round to value one the $B$ components of $\beta$ with the largest values. In practice, however, the $\beta$
returned by \texttt{Ipopt} was almost always already binary, so that the rounding step was rarely needed. Both solvers are run with default parameters and a time limit of
$1800$ seconds and, for a fair comparison, they are warm-started on the
 \lpnote{same starting value} used by our heuristics.
The two remaining baselines are the two variants of the recursive feature elimination (RFE) heuristic \lpnote{\gnote{introduced} in Section \ref{sec:related}.
Starting from the full set of $n$ features, RFE repeatedly eliminates a feature until the budget $B$ is reached on the basis of a specific ranking.
 The resulting binary vector $\beta$ is thus a feasible point
of~\eqref{prob:DFS-SVM}, satisfying $\sum_{j=1}^{n}\beta_j=B$.
The ranking procedure is carried out using the objective  $f(\alpha,\beta)$: \gnote{at each iteration
we discard the feature whose removal yields the smallest value of $f$.} \texttt{RFE1} and \texttt{RFE2} differ
 in how the candidate removals are
scored at each iteration.
\texttt{RFE1}, \gnote{the faster of the two}, keeps the dual variables $\alpha$ fixed at the solution of the
current subset while ranking all the candidate removals\gnote{, as in the kernel extension suggested
in \cite{guyon2002gene},} requiring a single SVM fit per iteration, whereas \texttt{RFE2} recomputes $\alpha$ for
every feature removal, being more accurate, but
requiring as many fits as active features at each iteration.   Neither variant requires a starting value $\beta^{\text{in}}$.
}
In Section~\ref{sec:certified}, we first compare this single-run version %
on instances for which  global optima can be certified. We then report the multistart results on much harder optimization instances for which optimality certification is, in general, computationally intractable.}

\subsubsection{\textcolor{revblue}{Instances with certified optima}}\label{sec:certified}
{\color{revblue}
We now evaluate the proposed heuristics in absolute terms, comparing their solutions against a
\emph{certified global optimum} of problem~\eqref{prob:DFS-SVM}.

The natural baseline for global optimality would be a
spatial branch-and-bound solver such as \texttt{Baron} \cite{tawarmalani2004global}. In a preliminary analysis we found it inadequate
to this aim, for two reasons. First, \texttt{Baron} cannot even return a feasible solution on any of the full instances, i.e., those with all the samples, within a two-hour time limit. Second, even on
reduced instances, with a  number of samples ranging from 10 to 200, its spatial branch-and-bound returns extremely large optimality gaps, up to the order of $10^{6}\%$, within the time limit
(see Appendix~\ref{sec:baron} for the detailed \texttt{Baron} comparison on reduced instances), which are therefore uninformative as certificates of global optimality. For these
reasons, we abandoned \texttt{Baron} as the reference for optimality and we resort to \gnote{a highly expensive enumerative approach}.

\lpnote{Indeed, for any given feature subset of dimension $B$, namely a fixed value $\beta$, the
problem \eqref{prob:DFS-SVM} reduces to a \emph{convex} SVM dual. We can therefore enumerate \emph{all}
$\binom{n}{B}$ cardinality-$B$ feature subsets, solve the corresponding convex dual exactly and take the
smallest objective value, which is the certified global optimum, obtained independently of any MINLP solver.}
This enumerative procedure is computationally burdensome, \lpnote{even though the subproblems are
mutually independent and can be solved in parallel: the enumeration remains affordable only up to a
dataset-dependent budget $B_{\max}$, reported in Table~\ref{tab:cert-coverage} in
Appendix~\ref{app:cert}, where we also motivate the restriction of the budget range to
$B\le\lfloor n/2\rfloor$.
It has been \gnote{run} on a single high-performance cloud computer (an AWS EC2
\texttt{c6i.32xlarge}: $128$ vCPUs over $64$ Intel Xeon Ice~Lake cores, $256$\,GB of RAM).}

 We highlight that the same workload on the $8$-core computer (Apple~M1) used for all the other experiments
would have been computationally intractable, in that it would have required several weeks of computation. Exact certification is therefore impractical as a general-purpose
tool, and we used it here only to prove the effectiveness of the fast heuristics proposed.

\begin{table}[!ht]
\centering\color{revblue}\renewcommand\arraystretch{1.2}\footnotesize
\resizebox{\textwidth}{!}{%
\begin{tabular}{lcc|cc|cc|cc|cc}
\toprule
& & & \multicolumn{2}{c|}{\texttt{Ipopt}} & \multicolumn{2}{c|}{\texttt{Knitro}} & \multicolumn{2}{c|}{\texttt{LS$^*$}} & \multicolumn{2}{c}{\AOsub} \\
dataset & $n$ & \# pbs & gap\,\% & \# solved & gap\,\% & \# solved & gap\,\% & \# solved & gap\,\% & \# solved \\
\midrule
Wholesale  & 7  &9 & 1.29  & 5/9  & 0.96  & 6/9  & \textbf{0.00}  & \textbf{9/9} & \textbf{0.00} & \textbf{9/9} \\
Diabetes   & 8  &12 & 2.26  & 0/12  & 2.78  & 1/12  & 0.02  & 11/12 & \textbf{0.00} & \textbf{12/12} \\
BCW        & 9  &12 & 12.99 & 3/12  & 11.21 & 3/12  & \textbf{0.00}  & \textbf{12/12} & \textbf{0.00} & \textbf{12/12} \\
Cleveland  & 13 &18 & 11.15 & 3/18  & 9.12  & 5/18  & 0.24  & 17/18 & \textbf{0.00} & \textbf{18/18} \\
Parkinsons & 22 &33 & 2.65  & 6/33  & 2.52  & 9/33  & 0.27  & 24/33 & \textbf{0.10} & \textbf{29/33} \\
German     & 24 &27 & 2.03  & 11/27 & 3.66  & 8/27  & 0.04  & 24/27 & \textbf{0.00} & \textbf{27/27} \\
BCD        & 30 &36 & 8.97  & 8/36  & 12.01 & 10/36 & 17.20 & 16/36 & \textbf{8.80} & \textbf{24/36} \\
Ionosphere & 33 &30 & 12.92 & 2/30  & 15.14 & 2/30  & 1.87  & 18/30 & \textbf{1.22} & \textbf{23/30} \\
Sonar      & 60 &18 & 7.56  & 1/18  & 9.91  & 0/18  & 0.95  & 11/18 & \textbf{0.19} & \textbf{15/18} \\
\midrule
Overall & &195 & 7.10 & 39/195 & 8.14 & 44/195 & 3.63 & 142/195 & \textbf{1.85} & \textbf{169/195} \\
\bottomrule
\end{tabular}}
\caption{
\lpnote{Polynomial kernels: comparison over the optimal certified
instances with  $B\le B_{\max}$.
``gap\,\%'' = mean gap to the certified optimum; ``\# solved'' = number of instances in which the method reaches
the optimum, over the certified instances. The ``Overall'' row
reports the average of the gaps and the number of instances solved to optimality over all instances. Bold marks the best value
in each row.}}
\label{tab:cert-poly}
\end{table}

\paragraph{Results.}
For each dataset, we construct MINLP instances \eqref{prob:DFS-SVM} with  $C=10$, $\gamma=0.1$ and
\begin{itemize}
\item polynomial kernels of degree $d=2,3,5$ with $c=1$, and the Gaussian kernel;
\item  budget $B\in [1, B_{\max}]$
\end{itemize}
for a total of 195 instances for the polynomial kernels, and 63 instances for the Gaussian kernel.

We run the two heuristics \texttt{LS$^*$} and \AOsub and  the two general-purpose solvers \texttt{Ipopt} and \texttt{Knitro} on the set of instances solved to certified optimality by the complete enumeration procedure.
For each pair (dataset, algorithm) we report in
Table~\ref{tab:cert-poly} (polynomial kernels case) and Table~\ref{tab:cert-rbf} (Gaussian
kernel case):
\begin{itemize}
\item the overall number of instances \gnote{(\# pbs);}
    \item the mean gap to the certified optimum \gnote{(avg gap\%);}
    \item the number of instances \gnote{on which the method reaches that optimum (\# solved).}
\end{itemize}

\begin{table}[!ht]
\centering\color{revblue}\renewcommand\arraystretch{1.2}\footnotesize
\begin{tabular}{lcc|cc|cc|cc}
\toprule
 & & & \multicolumn{2}{c|}{\texttt{Ipopt}} & \multicolumn{2}{c|}{\texttt{Knitro}} & \multicolumn{2}{c}{\texttt{LS$^*$}} \\
dataset & $n$ & \# pbs & gap\,\% & \# solved & gap\,\% & \# solved & gap\,\% & \# solved \\
\midrule
Wholesale  & 7  & 3  & 1.64  & 1/3  & 1.14  & 2/3  & \textbf{0.00} & \textbf{3/3}  \\
Diabetes   & 8  & 4  & 0.82  & 1/4  & 3.20  & 0/4  & \textbf{0.00} & \textbf{4/4}  \\
BCW        & 9  & 4  & 7.71  & 1/4  & 4.79  & 1/4  & \textbf{0.00} & \textbf{4/4}  \\
Cleveland  & 13 & 6  & 10.03 & 1/6  & 10.03 & 1/6  & \textbf{0.00} & \textbf{6/6}  \\
Parkinsons & 22 & 11 & 1.41  & 2/11 & 2.12  & 3/11 & \textbf{1.16} & \textbf{6/11} \\
German     & 24 & 9  & 1.68  & 3/9  & 2.19  & 2/9  & \textbf{0.03} & \textbf{7/9}  \\
BCD        & 30 & 11 & \textbf{5.96}  & 2/11 & 11.27 & 2/11 & 5.98 & \textbf{5/11} \\
Ionosphere & 33 & 9  & 16.73 & 0/9  & 13.19 & 0/9  & \textbf{3.65} & \textbf{5/9}  \\
Sonar      & 60 & 6  & 4.66  & 0/6  & 9.40  & 0/6  & \textbf{2.20} & \textbf{3/6}  \\
\midrule
Overall & & 63 & 5.93 & 11/63 & 6.95 & 11/63 & \textbf{1.98} & \textbf{43/63} \\
\bottomrule
\end{tabular}
\caption{\gnote{Gaussian (RBF) kernel: comparison against the certified global optimum, over the certified
instances with $B\le B_{\max}$.
\AOsub is not applicable (it is restricted to polynomial kernels), so only \texttt{LS$^*$}
and the general-purpose solvers \texttt{Ipopt} and \texttt{Knitro} are reported. Columns as in Table~\ref{tab:cert-poly}.}}
\label{tab:cert-rbf}
\end{table}

\begin{table}[!ht]
\centering\color{revblue}\renewcommand\arraystretch{1.2}\footnotesize
\begin{tabular}{lc|cc|cc|cc|cc}
\toprule
 & & \multicolumn{4}{c|}{polynomial kernels} & \multicolumn{4}{c}{Gaussian kernel} \\
 & & \multicolumn{2}{c|}{\texttt{RFE1}} & \multicolumn{2}{c|}{\texttt{RFE2}} & \multicolumn{2}{c|}{\texttt{RFE1}} & \multicolumn{2}{c}{\texttt{RFE2}} \\
dataset & $n$ & gap\,\% & \# solved & gap\,\% & \# solved & gap\,\% & \# solved & gap\,\% & \# solved \\
\midrule
Wholesale  & 7  & 0.81  & 7/9   &  \textbf{0.00}  &  \textbf{9/9}   &  \textbf{0.00}  &  \textbf{3/3}  &  \textbf{0.00} &  \textbf{3/3} \\
Diabetes   & 8  & 3.97  & 2/12  & 0.05  & 10/12 & 0.33  & 2/4  & 0.07 & 3/4 \\
BCW        & 9  & 20.44 & 3/12  & 6.78  & 5/12  & 12.24 & 1/4  & 2.09 & 3/4 \\
Cleveland  & 13 & 18.96 & 3/18  & 1.90  & 12/18 & 20.84 & 1/6  & 0.57 & 4/6 \\
Parkinsons & 22 & 5.34  & 5/33  & 1.83  & 12/33 & 1.61  & 3/11 & \textbf{1.11} & 3/11 \\
German     & 24 & 0.89  & 17/27 & 0.55  & 20/27 & 0.80  & 5/9  & 0.51 & 6/9 \\
BCD        & 30 & 45.81 & 3/36  & 55.63 & 0/36  & 37.03 & 0/11 & 57.45 & 0/11 \\
Ionosphere & 33 & 23.17 & 0/30  & 14.64 & 2/30  & 19.99 & 0/9  & 8.58 & 0/9 \\
Sonar      & 60 & 15.29 & 0/18  & 7.63  & 1/18  & 7.55  & 1/6  & 4.90 & 0/6 \\
\midrule
Overall & & 17.75 & 40/195 & 14.21 & 71/195 & 13.22 & 16/63 & 12.18 & 22/63 \\
\bottomrule
\end{tabular}
\caption{\gnote{Recursive-feature-elimination baselines on the same certified instances of
Tables~\ref{tab:cert-poly}--\ref{tab:cert-rbf}, again with $B\le B_{\max}$.
Columns as in
Table~\ref{tab:cert-poly}. Here bold marks the best value among all the methods
reported in Tables~\ref{tab:cert-poly}--\ref{tab:cert-rfe} for the corresponding kernel family,
and not the best of the two baselines in each row.}}
\label{tab:cert-rfe}
\end{table}

\begin{figure}[!ht]\centering\color{revblue}
  \includegraphics[width=\linewidth]{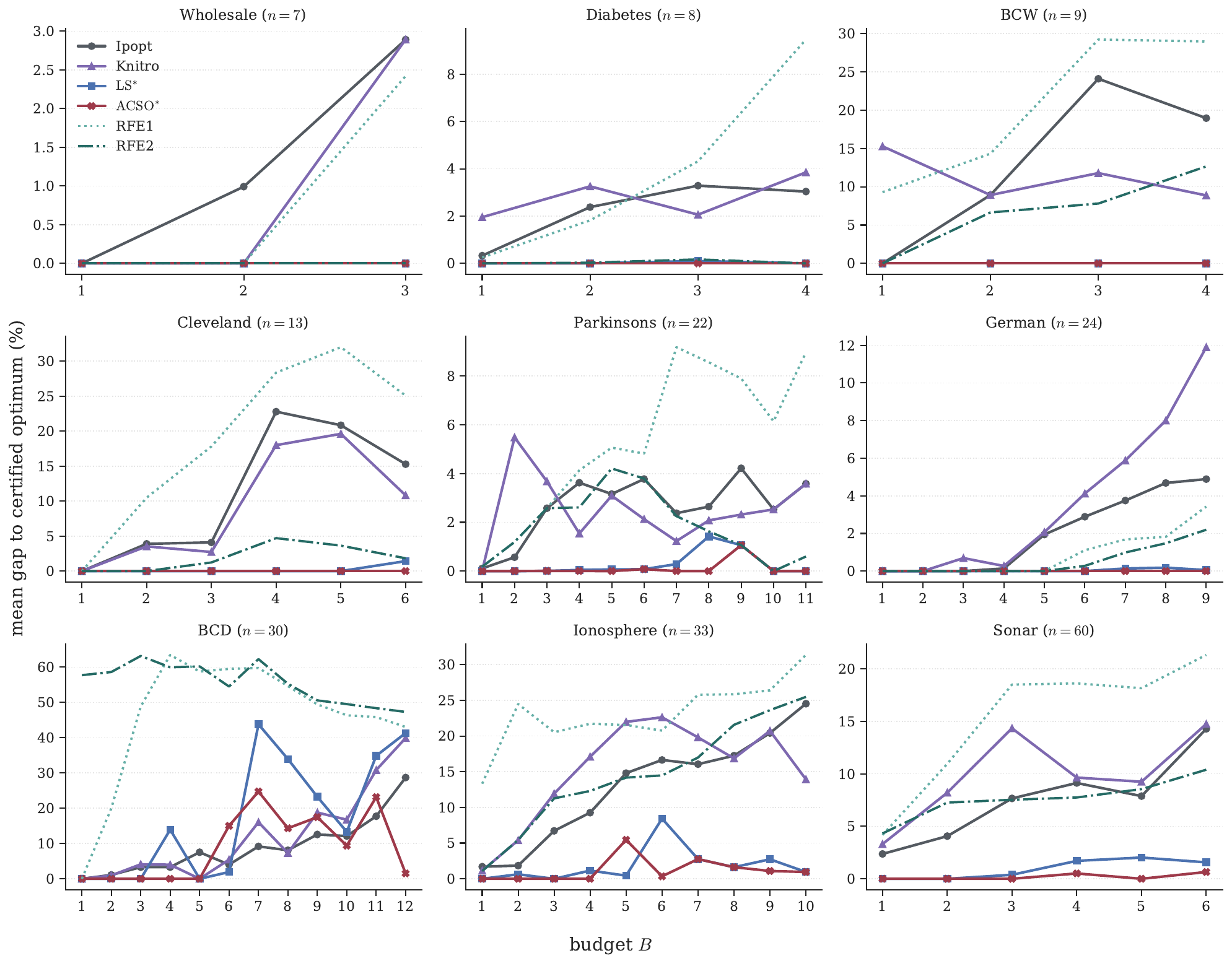}
  \caption{\gnote{Mean gap over the three polynomial kernels to the certified global optimum as the cardinality budget $B$ grows, for each dataset. The budget ranges from $1$ to the certified value $B_{\max}$ reported in Table~\ref{tab:cert-coverage}.}
}
  \label{fig:cert-gapB}
\end{figure}

\begin{figure}[!ht]\centering\color{revblue}
  \includegraphics[width=\linewidth]{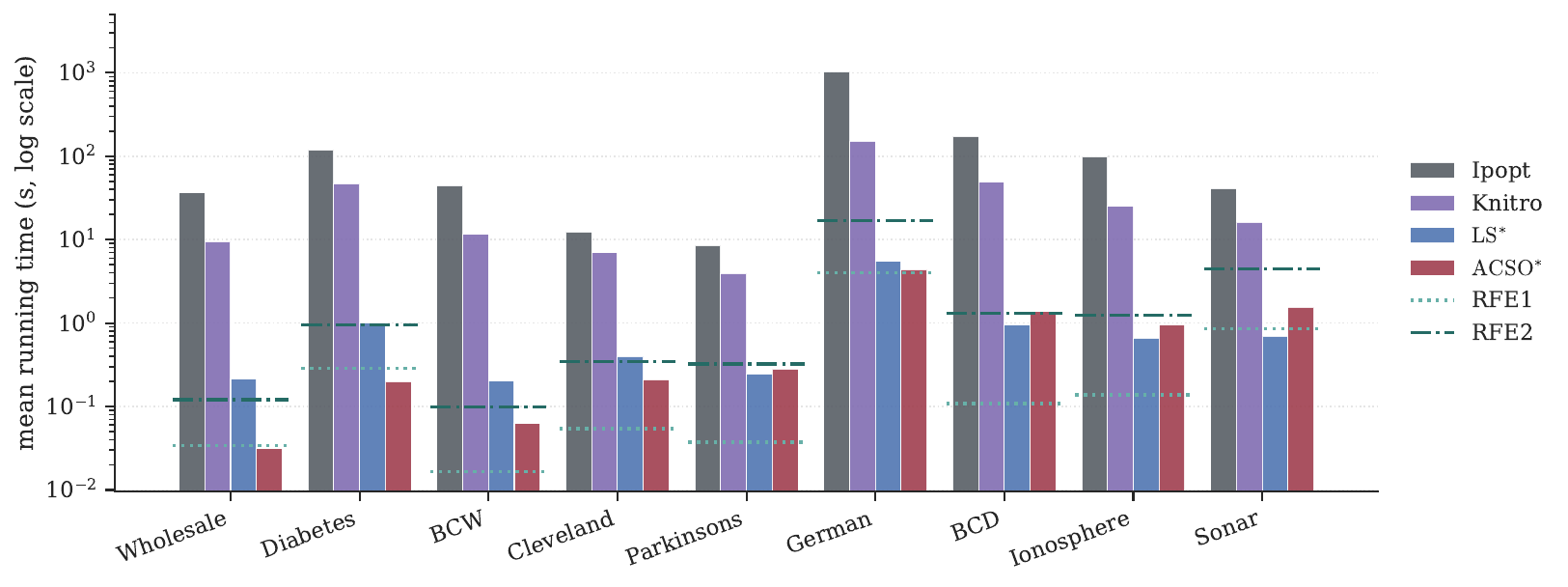}
  \caption{\gnote{
  Mean running time  (logarithmic
  scale) for each dataset on the certified instances with polynomial kernels; for each
  dataset the mean is taken over its budgets $B=1,\dots,B_{\max}$ and the three polynomial degrees. The dotted lines mark the
  \texttt{RFE1} and \texttt{RFE2} heuristics.}  }
  \label{fig:cert-time}
\end{figure}
}

\gnote{On the polynomial kernels the two
heuristics are far closer to the certified optimum than the solvers: \AOsub attains a
mean gap of $1.85\%$ and reaches the certified optimum in $169/195$ instances, and \texttt{LS$^*$} $3.63\%$
($142/195$), whereas \texttt{Knitro} and \texttt{Ipopt} stop at $8.14\%$ and $7.10\%$ and reach the optimum in only $44/195$ and $39/195$, respectively. The same
holds on the Gaussian kernel, where the decomposition does not apply and \texttt{LS$^*$} ($1.98\%$,
$43/63$) clearly outperforms both \texttt{Knitro} ($6.95\%$, $11/63$) and \texttt{Ipopt} ($5.93\%$, $11/63$).
For completeness, Table~\ref{tab:cert-rfe} evaluates on the same instances the two
recursive-feature-elimination baselines \texttt{RFE1} and \texttt{RFE2}: although they are as fast as our heuristics
(mean selection time below $3$\,s), their solutions are far from the certified optimum, with mean gaps of
$17.8\%$ and $14.2\%$ on the polynomial kernels ($13.2\%$ and $12.2\%$ on the Gaussian one), even farther
than the general-purpose solvers. Figure~\ref{fig:cert-gapB} shows, per dataset, how this gap evolves as the budget $B$ grows: on
most datasets the heuristics reach the certified optimum over the whole budget range, and only a few
hard instances leave a gap larger than a fraction of a percent}.

\gnote{Beyond solution quality, the two heuristics are also one to two orders of magnitude faster than the solvers (Figure~\ref{fig:cert-time}): on the polynomial instances the mean running time is $1.3$\,s for \texttt{LS$^*$} and $1.2$\,s for \AOsub{}, against $62$\,s for \texttt{Knitro} and $239$\,s for \texttt{Ipopt}; on the Gaussian kernel \texttt{LS$^*$} takes $2.1$\,s, against $37$\,s and $111$\,s for \texttt{Knitro} and \texttt{Ipopt}. In Figure~\ref{fig:cert-time} the running times are reported per dataset on a \emph{logarithmic} time axis (each gridline is a ten-fold increase): this is needed because the times span several orders of magnitude, from a fraction of a second for the heuristics to hundreds of seconds for the solvers, which a linear axis could not display legibly.}

\subsubsection{\textcolor{revblue}{Multistart experiments}}\label{sec:noncert}

\begin{figure}[!ht]\centering
  \includegraphics[width=\linewidth]{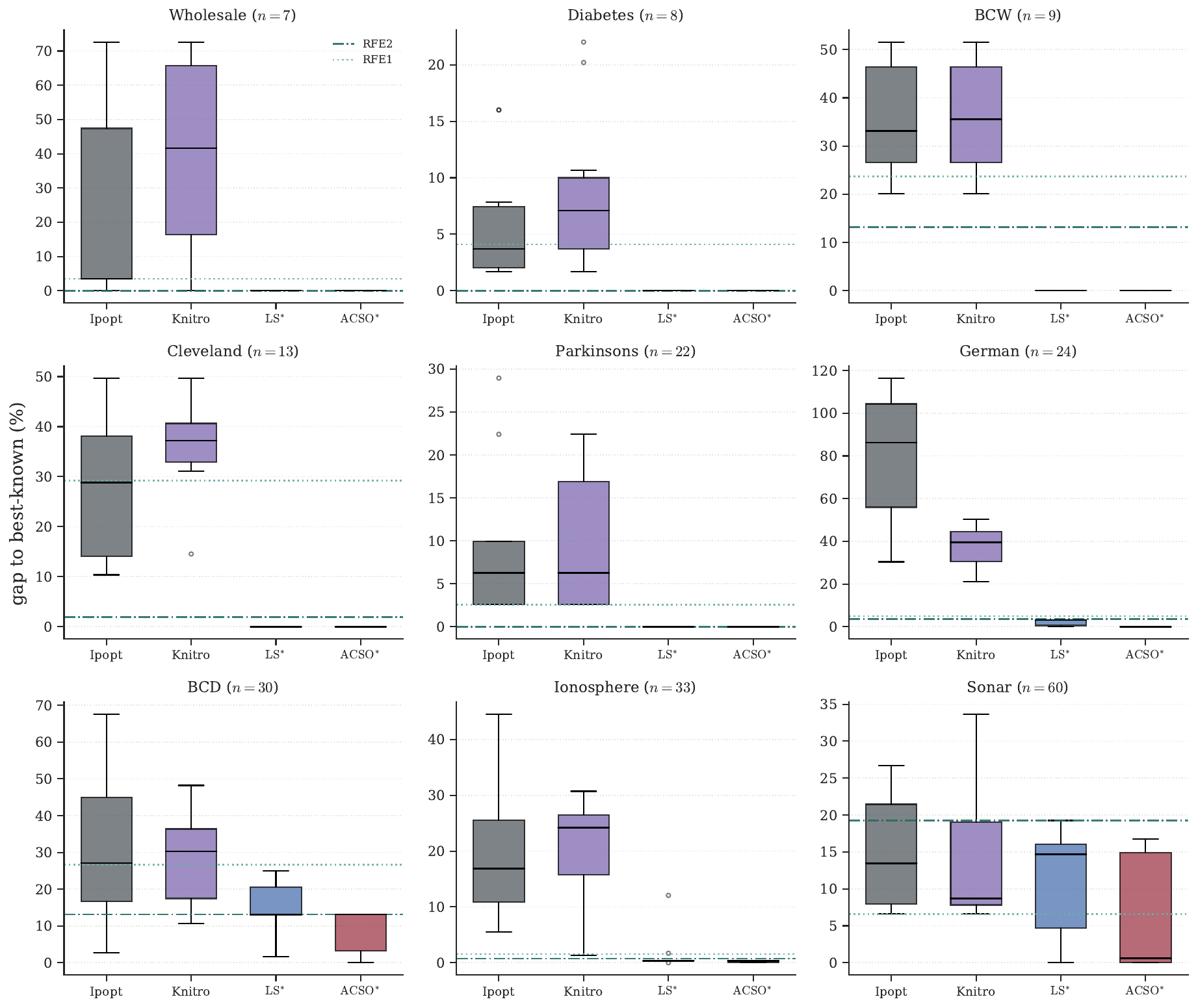}
  \caption{\textcolor{revblue}{Box plot of the objective gap to the best-known value across the $10$  random starts;
  the dotted lines mark the deterministic \texttt{RFE2}/\texttt{RFE1}. \textcolor{revblue}{For the first five
datasets ($n\le 22$) this gap equals the certified optimality gap, because the best-known objective function value coincides with the optimal one derived from the certification procedure of Section~\ref{sec:certified}.}}}
  \label{fig:b2-box}
\end{figure}

\begin{figure}[!ht]\centering
  \includegraphics[width=\linewidth]{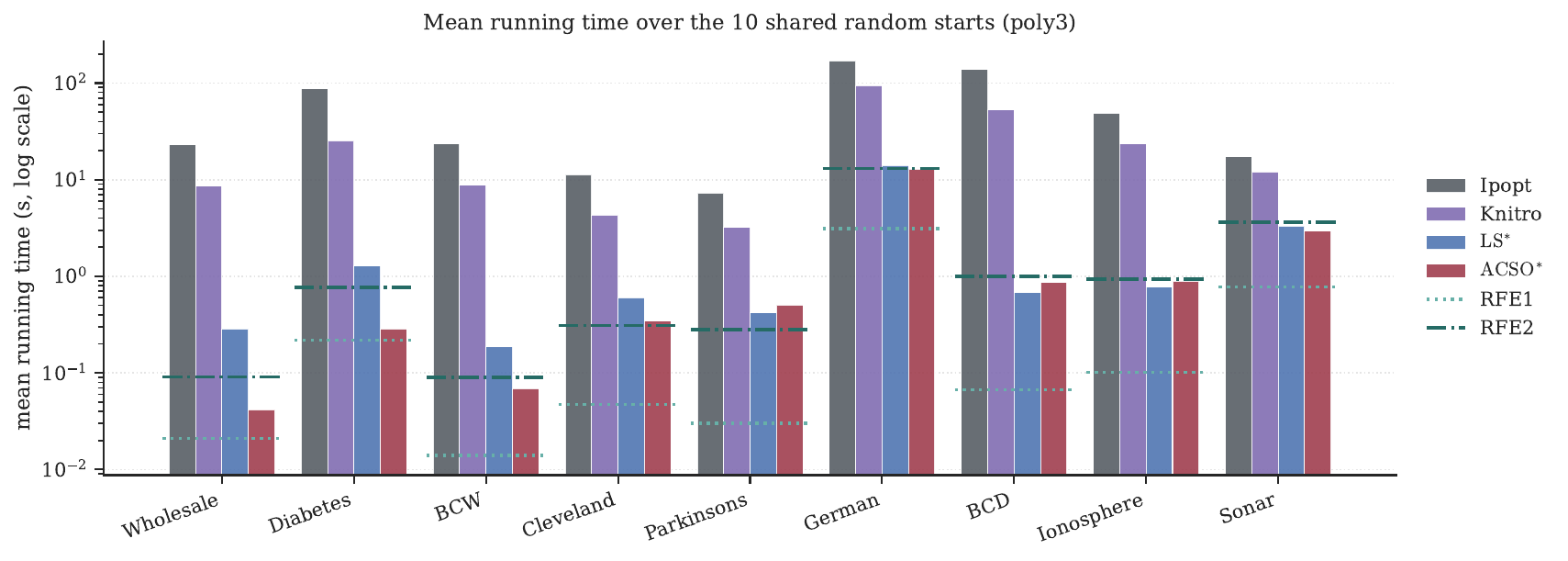}
  \caption{\textcolor{revblue}{Mean running time (\emph{logarithmic} scale) over the 10 runs on the $B=\lfloor n/2\rfloor$ instances  (cubic polynomial
  kernel).%
  The dotted lines mark
  the deterministic \texttt{RFE1}, \texttt{RFE2}.}}
  \label{fig:b2-time}
\end{figure}

\begin{table}[!ht]
\centering\color{revblue}\renewcommand\arraystretch{1.2}
\begin{tabular}{lcccc}
\toprule
method & \# instances & mean gap to best-known (\%) & \# best-known reached & mean spread (\%) \\
\midrule
\texttt{Ipopt}  & 36 & 10.35 & 4/36  & 38.6 \\
\texttt{Knitro} & 36 &  9.88 & 7/36  & 34.1 \\
\texttt{RFE1}   & 36 & 13.99 & 3/36  & --   \\
\texttt{RFE2}   & 36 &  6.23 & 12/36 & --   \\
\texttt{LS$^*$} & 36 &  0.51 & 32/36 &  6.3 \\
\AOsub          & 27 & \textbf{0.21} & \textbf{26/27} & \textbf{4.4} \\
\bottomrule
\end{tabular}
\caption{\gnote{Multistart experiments at $B=\lfloor n/2\rfloor$, over $10$ shared random starting
points. ``mean gap'' $=$ gap to the best-known objective of the best of the $10$ starts, averaged over
the instances; ``\# best-known reached'' $=$ number of instances on which the best-known objective is
attained; ``mean spread'' $=$
average spread of the gap across the $10$ starts, defined in the text. \AOsub is reported over the $27$
polynomial instances, the only ones to which it applies; \texttt{RFE1} and \texttt{RFE2} are
deterministic and therefore have no spread.}}
\label{tab:b2-summary}
\end{table}

\begin{figure}[!ht]\centering
  \includegraphics[width=\linewidth]{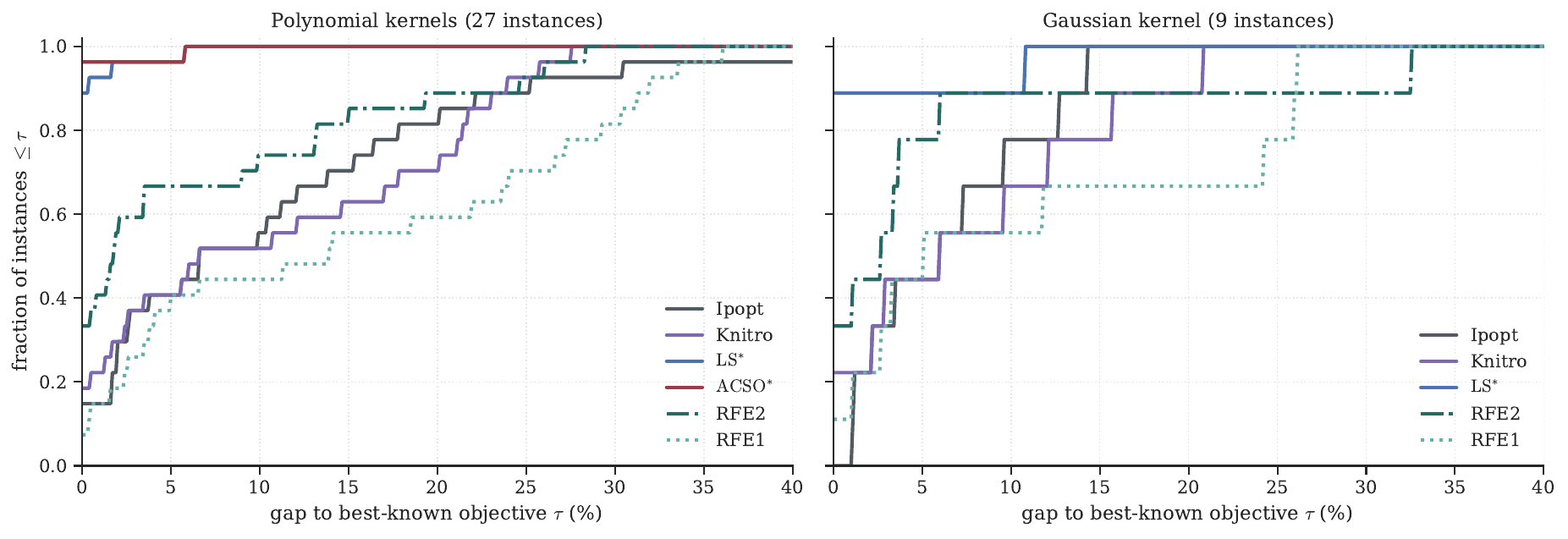}
  \caption{\textcolor{revblue}{Performance profiles on the $B=\lfloor n/2\rfloor$ instances (polynomial kernels,
  left; Gaussian kernel, right). For each method we plot the fraction of instances whose gap to the
  best-known objective does not exceed a threshold $\tau$ (\%); each method is represented by its best
  solution over the $10$ shared random starts. A curve rising steeply near $\tau=0$ denotes a method that is
  at or near the best on most instances. \AOsub is polynomial-only and hence absent from the Gaussian plot.}}
  \label{fig:b2-pp}
\end{figure}

{\color{revblue}

\lpnote{In Section~\ref{sec:certified} all methods \gnote{start} from the same deterministic
point and the budget varies in the interval $[1, B_{\max}]$.}
 We now fix the budget at $B=\lfloor n/2\rfloor$, and we analyse the dependence of
each method on the starting point.
 \lpnote{The value $B=\lfloor n/2\rfloor$ is }
 the combinatorially
hardest value, for which $\binom{n}{B}$ is largest (e.g.\ $\binom{60}{30}\approx 1.2\times 10^{17}$ for
Sonar) and a certified optimum is therefore not necessarily available\gnote{.} We compare our two heuristics \texttt{LS$^*$} and \AOsub\  against \texttt{Knitro}, \texttt{Ipopt},
\texttt{RFE1} and \texttt{RFE2}.

Every start-dependent method is run from \emph{the same} $10$ random feature
subsets $\beta^{\text{in}}$, so that no method benefits from a tailored initialization; %
\texttt{RFE1}/\texttt{RFE2} are deterministic and need no starting solution.
Whether a certified optimum is available at $B=\lfloor n/2\rfloor$ depends on the dataset. For the five
datasets with $n\le 22$ (Wholesale, Diabetes, BCW, Cleveland and Parkinsons) the
certified range of Section~\ref{sec:certified} reaches exactly $B_{\max}=\lfloor n/2\rfloor$
(Table~\ref{tab:cert-coverage}), so this budget is certified; for the remaining four (German, BCD, Ionosphere and Sonar, all with $n\ge 24$) the enumeration stops earlier,
$B_{\max}<\lfloor n/2\rfloor$, and no certified optimum is available at this budget. Each method is
therefore scored against the \emph{best-known} objective of every instance, i.e., the best value found by
any method from any start. On the five datasets with $n\le 22$, the best-known value
coincides \emph{exactly} with the certified optimum, and both \texttt{LS$^*$} and \AOsub attain it on all
of them, an independent confirmation that our heuristics reach the true optimum.

\gnote{Figure~\ref{fig:b2-box} reports the box plots of the objective-function gaps for all the methods, \texttt{Knitro}, \texttt{Ipopt}, \texttt{LS$^*$}, \AOsub, \texttt{RFE1} and \texttt{RFE2}, on the polynomial kernel with $d=3$.}
\texttt{LS$^*$} and \AOsub sit at or near the best-known value with a very small
spread across the $10$ starts, while \texttt{Knitro} and \texttt{Ipopt} leave large gaps and, crucially, are
\emph{strongly dependent on the starting point} (wide boxes). \texttt{RFE1}/\texttt{RFE2}, being
deterministic, are single points and remain far from the best.
In Figure~\ref{fig:b2-time}, running times are compared {for the instances with polynomial kernels of degrees $d=2,3,5$}.
The two heuristics are one to two orders of magnitude faster than the solvers (mean $2.1$--$2.9$\,s per run
against $25.7$\,s for \texttt{Knitro} and $54.0$\,s for \texttt{Ipopt}).

\gnote{Table~\ref{tab:b2-summary} summarizes all the multistart experiments over the whole test bed, which consists of
$36$ instances, one for each of the $9$ datasets combined with each of the $4$ kernels, all taken at
$B=\lfloor n/2\rfloor$. Of these, the $27$ instances built on a polynomial kernel are those on which
\AOsub is applicable, whereas all
the remaining methods are run on all $36$. For every instance a method is credited with reaching the
best-known objective when it matches the best value found by any method from any start. Our heuristics  \AOsub and \texttt{LS$^*$}  reach the best-known objective on almost every instance, while the two
general-purpose solvers do so only occasionally and the \texttt{RFE} baselines lie in between. Denoting
by $f^{\mathrm{bk}}$ the best-known objective of an instance and by $\mathrm{gap}(f)=100\,
(f-f^{\mathrm{bk}})/|f^{\mathrm{bk}}|$ the gap of a solution of value $f$ to the best-known, the last column quantifies
the dependence on the starting point as the average, over the instances, of
$\mathrm{gap}(f_{\mathrm{worst}})-\mathrm{gap}(f_{\mathrm{best}})$, where $f_{\mathrm{best}}$ and
$f_{\mathrm{worst}}$ are the best and the worst of the $10$ starts. This spread is of a few percentage
points for our heuristics and of several tens for the solvers, which is the same picture conveyed by
the box plots of Figure~\ref{fig:b2-box}. As shown in Figure~\ref{fig:b2-time}, the heuristics are also
one to two orders of magnitude faster.}

Figure~\ref{fig:b2-pp} aggregates all instances through performance profiles in the sense
of~\cite{dolan2002benchmarking}: the fraction of instances solved within a relative gap $\tau$ to the
best-known objective. On both the polynomial and the Gaussian kernels, \AOsub and \texttt{LS$^*$} reach the
best-known value on almost all instances already at $\tau\!\approx\!0$, whereas \texttt{Knitro},
\texttt{Ipopt} and the \texttt{RFE} baselines need much larger thresholds to catch up.}

\section{Conclusion}
\label{sec: 6}

We have addressed the challenging problem of embedding feature selection in nonlinear SVMs under a strict cardinality constraint, formulating it as a Mixed-Integer Nonlinear Programming problem. While the literature on feature selection for linear SVMs is extensive, mixed-integer approaches for nonlinear SVMs are rare due to their substantially greater computational complexity. Our work fills this gap by introducing two novel algorithmic frameworks.

First, we proposed a local search metaheuristic applicable to general nonlinear kernels, capable of efficiently exploring the binary search space and delivering high-quality solutions. Second, we developed a decomposition framework based on alternating optimization between continuous and binary variables; for polynomial kernels, we exploited a submodular reformulation of the binary subproblem, enabling the use of scalable submodular maximization algorithms. To the best of our knowledge, this is the first application of submodular optimization techniques to interpretable nonlinear SVM models.

Extensive numerical experiments demonstrate that our methods consistently outperform state-of-the-art solvers and heuristic baselines, achieving superior solutions within substantially reduced computational times. These results open a promising direction for efficiently training  nonlinear SVMs while guaranteeing exact control over sparsity.

Future work will focus on extending the decomposition strategy to Gaussian kernels, adapting the proposed frameworks to other machine learning tasks such as regression and clustering, and integrating them into the model of~\cite{MARGOT} to enable optimal decision trees with nonlinear SVM-based splits.

\lpnote{\section*{Data availability and reproducibility}
The results presented in this paper are fully reproducible and the code, together with the datasets used in this study, is available in the public repository at \texttt{https://github.com/f-donofrio/fs-nl-svm}.
}

\printbibliography

\appendix

\section{Appendix}

\subsection{Proofs of Proposition~\ref{prop: submodularity and monotonicity of F} and Proposition~\ref{prop: H sub}}\label{app:proofs}
We prove in the following Propositions  \ref{prop: submodularity and monotonicity of F} and \ref{prop: H sub}.
Preliminarily, we make the following observation which allows us to prove \textcolor{revblue}{Proposition}~\ref{prop: submodularity and monotonicity of F}.

\begin{observation}
\label{observation: c poly kernel = 0 poly kernel}

{Given a dataset with features vectors $x^i \in \mathbb{R}^{n}, i \in [m]$ } and a subset of features $S \in N$, a polynomial kernel matrix $K_S$, with components $K_{ih} = k(x^i_S, x^h_S) = \left (c + \sum_{j \in S} \gamma x^i_j x^h_j \right)^d$, $i,h\in [m]$, $c \neq 0$, can be expressed as a polynomial kernel $\Tilde{K}(\Tilde{x}, \Tilde{x})$ with components $\Tilde{k}(\Tilde{x}^i_S, \Tilde{x}^h_S ) = \left ( \sum_{j \in S} \gamma \Tilde{x}^i_j \Tilde{x}^h_j \right)^d$, $i,h\in [m]$ where $\Tilde{x}^i\in R^{(n+1)}, i\in[m],$ is a new  {features vector} with values $\Tilde{x}^i_j$, $j\in   [n]\cup\{0\}$,  such that $\Tilde{x}^i_0 = \sqrt{\frac c \gamma} $ and $\Tilde{x}^i_{j} = x^i_{j}$, for $j \in [n]$. Indeed we have:
\begin{equation*}
k(x^i_S, x^h_S) = \left (c + \sum_{j \in S} \gamma x^i_j x^h_j \right)^d =  \left (\sum_{j \in S \cup \{0\}} \gamma \Tilde{x}^i_j \Tilde{x}^h_j \right )^d = \Tilde{k}(\Tilde{x}^i_S, \Tilde{x}^h_S).
\end{equation*}

Let us now denote by $\Tilde F$ the same function $F$ but defined over  {these new features vectors}, i.e., over sets $\Tilde{S} \in \Tilde{N}$, where $\Tilde{N}= 2^{[n]\cup\{0\}}$. We write  $F_{|_{c=0}}$  or $\Tilde F_{|_{c=0}}$ to indicate that the hyperparameter $c$ of the polynomial is set to $0$. Thus, we have that for any set $S \in N$:
\begin{equation*}
F(S) = \sum_{i=1}^m\sum_{h=1}^m \bar\alpha_i\bar \alpha_h y^iy^h \left (c + \sum_{j \in S} \gamma x^i_j x^h_j \right )^d  =  \sum_{i=1}^m\sum_{h=1}^m \bar\alpha_i\bar \alpha_h y^iy^h \left (\sum_{j \in S \cup \{0\}} \gamma \Tilde{x}^i_j \Tilde{x}^h_j \right )^d = \Tilde{F} (S \cup \{0\})_{|_{c=0}}.
\end{equation*}

\end{observation}

In other words, by Observation \ref{observation: c poly kernel = 0 poly kernel}, any function $F$ with the $c$ hyperparameter different from $0$, can be redefined as a function $\Tilde{F}$ with $c$ set to $0$, by just modifying the input dataset.

\textbf{Proof of Proposition \ref{prop: submodularity and monotonicity of F}.}

\begin{proof}
Firstly we prove supermodularity. Using Observation \ref{observation: c poly kernel = 0 poly kernel}, we can focus on the case in which the coefficient $c$ of the polynomial kernel is $0$. Thus the set function of interest is
$$
F(S)_{\mid c=0}=\sum_{i=1}^m\sum_{h=1}^m\bar\alpha_i\bar\alpha_h y^iy^h\Bigl(\sum_{j\in S}\gamma x^i_j x^h_j\Bigr)^d
=\sum_{i=1}^m\sum_{h=1}^m a_{ih}\Bigl(\sum_{j\in S}g^{ih}_j\Bigr)^d,
$$
where we write $a_{ih}:=\bar\alpha_i\bar\alpha_h y^iy^h$ and $g^{ih}_j:=\gamma x^i_j x^h_j$. We must show that, for any $A\subseteq B$ and any feature $e\in[n]\setminus B$,
$$
F(A\cup\{e\})_{\mid c=0}-F(A)_{\mid c=0}\;\le\;F(B\cup\{e\})_{\mid c=0}-F(B)_{\mid c=0},
$$
which is equivalent to
$$
F(B\cup\{e\})_{\mid c=0}-F(B)_{\mid c=0}-\bigl(F(A\cup\{e\})_{\mid c=0}-F(A)_{\mid c=0}\bigr)\ge0.
$$

Define the marginal gain $\Delta_e(S):=F(S\cup\{e\})_{\mid c=0}-F(S)_{\mid c=0}$ for $e\notin S$. Using the multinomial expansion
\begin{equation}\label{eq:multinomial-expansion}
\Bigl(\sum_{j\in S}g^{ih}_j\Bigr)^d=\sum_{j_1,\dots,j_d\in S}\prod_{\ell=1}^d g^{ih}_{j_\ell},
\end{equation}

a direct binomial expansion yields

 \begin{align}
\Delta_e(S) &  = \ds \sum_{i=1}^m\sum_{h=1}^m a_{ih} \left ( \sum_{j_1, \dots, j_d \in S \cup \{e\}} \prod_{\ell = 1}^d g_{j_\ell}^{ih} - \sum_{j_1, \dots, j_d \in S } \prod_{\ell = 1}^d g_{j_\ell}^{ih}\right )  \nonumber\\
 & = \ds \sum_{i=1}^m\sum_{h=1}^m a_{ih} \left ( \sum_{\substack{j_1, \dots, j_d \in S \cup \{e\} \nonumber\\ \text{at least one $j_{\ell} = e$}} } \prod_{\ell = 1}^d g_{j_\ell}^{ih}  \right ) \\
 & = \sum_{i=1}^m\sum_{h=1}^m a_{ih} \left (\sum_{k=1}^{d}\binom{d}{k}\bigl(g^{ih}_e\bigr)^k\!\!{\sum_{j_1,\dots,j_{d-k}\in S}}\prod_{\ell=1}^{d-k}g^{ih}_{j_\ell} \right ).\label{eq:increase}
 \end{align}
To see why the last equation holds, observe that each selection of $j_1,\dots, j_d \in S\textcolor{revblue}{\cup\{e\}}$ so that at least one $j_\ell$ is equal to $e$ is univocally identified by: the set of $k \geq 1$ indices $i_1,\dots,i_k$ from $\{1,\dots, d\}$  such that $j_{i_s} \textcolor{revblue}{=e}$ for $s=1,\dots, k$; the ordered selection (with repetitions) of the remaining $d-k$ entries from $S$.
We thus have that

 \begin{align}
  \Delta_e(B) - \Delta_e(A)  & =  \sum_{i=1}^m\sum_{h=1}^m a_{ih}\sum_{k=1}^{d}\binom{d}{k}\bigl(g^{ih}_e\bigr)^k\!\!\left(\sum_{{ j_1,\dots,j_{d-k} \in B}}\prod_{\ell=1}^{d-k}g^{ih}_{j_\ell}\right ) - \sum_{i=1}^m\sum_{h=1}^m a_{ih}\sum_{k=1}^{d}\binom{d}{k}\bigl(g^{ih}_e\bigr)^k\!\!\left(\sum_{{ j_1,\dots,j_{d-k} \in A} }\prod_{\ell=1}^{d-k}g^{ih}_{j_\ell}\right) \nonumber  \\
    & =  \sum_{i=1}^m\sum_{h=1}^m a_{ih}\sum_{k=1}^{d}\binom{d}{k}\bigl(g^{ih}_e\bigr)^k\!\!\left(\sum_{{ j_1,\dots,j_{d-k} \in B}}\prod_{\ell=1}^{d-k}g^{ih}_{j_\ell} - \sum_{{ j_1,\dots,j_{d-k} \in A} }\prod_{\ell=1}^{d-k}g^{ih}_{j_\ell}\right)  \nonumber \\
   & =  \sum_{i=1}^m\sum_{h=1}^m a_{ih}\sum_{k=1}^{d}\binom{d}{k}\bigl(g^{ih}_e\bigr)^k\!\!\left(\sum_{\substack{{j_1,\dots,j_{d-k}\in B} \\
 \text{at least one }j_\ell\in B\setminus A}}\prod_{\ell=1}^{d-k}g^{ih}_{j_\ell}\right) \nonumber\\
   &  = \sum_{i=1}^m\sum_{h=1}^ma_{ih} \left[
  \sum_{k=1}^{d}\binom{d}{k}\left(
  \sum_{t=1}^{d-k}{d-k \choose t} \bigl(g^{ih}_e\bigr)^{k} \sum_{j_1,\dots,j_t \in B \setminus A}  \prod_{\ell=1}^t g_{j_\ell}^{ih} \sum_{j_{t+1},\dots,j_{d-k}\in A} \prod_{q=t+1}^{d-k} g_q^{ih} \right) \right] \nonumber \\
&  =
  \sum_{i=1}^m\sum_{h=1}^ma_{ih} \left[ \sum_{k=1}^{d}\binom{d}{k}\left(
  \sum_{t=1}^{d-k} {d-k \choose t} \underbrace{\bigl(g^{ih}_e\bigr)^{k}}_{k_1(x^i,x^h)} \underbrace{(\sum_{j\in B \setminus A}g^{ih}_j\Bigr)^{t}}_{k_2(x^i,x^h)} \underbrace{(\sum_{j\in A}g^{ih}_j\Bigr)^{d-k-t}}_{k_3(x^i,x^h)} \right) \right]. \label{eq:last}
\end{align}
To see why the second-to-last equation above holds, observe that each selection of $j_1,\dots, j_d \in B$ so that at least one $j_\ell$ belongs to $B\setminus A$ is univocally identified by: the number $t\geq 1$ of entries $j \in \{j_1,\dots,j_{d-k}\}$ so that $j \in B\setminus A$; \textcolor{revblue}{which $t$ of the $d-k$ positions $\{j_1,\dots,j_{d-k}\}$ take their value in $B\setminus A$ (there are exactly $\binom{d-k}{t}$ such choices, which is the origin of the binomial factor $\binom{d-k}{t}$ in~\eqref{eq:last}), the remaining $d-k-t$ positions taking their value in $A$;} the ordered selection (with repetitions) of \textcolor{revblue}{these $t$ entries from $B\setminus A$ and of} the remaining $d-k-t$ entries from $A$. \textcolor{revblue}{For a fixed such choice, summing over all ordered selections with repetition yields the factors $(\sum_{j\in B\setminus A}g^{ih}_j)^{t}$ and $(\sum_{j\in A}g^{ih}_j)^{d-k-t}$} by~\eqref{eq:multinomial-expansion}.
Now observe that $k_1(\cdot,\cdot),k_2(\cdot,\cdot),k_3(\cdot,\cdot)$ are kernels.  Using the following proposition, we deduce that the term within the square bracket of~\eqref{eq:last} is a kernel $k_4(\cdot, \cdot)$.

\begin{proposition}[\cite{shawe2004kernel}]
Given two kernel functions, $k_1(\cdot, \cdot)$ and $k_2(\cdot, \cdot)$, and two nonnegative values $v_1, v_2 \in \R^+$, the sum $v_1k_1(\cdot, \cdot) + v_2 k_2(\cdot, \cdot)$ and the product $v_1 k_1(\cdot,\cdot) \cdot v_2 k_2(\cdot, \cdot)$ , are also kernel functions.
\label{prop: product of kernels}
\end{proposition}

Thus, we can use the definition of $a_{ih}$ and write
\begin{align*}
 \Delta_e(B) - \Delta_e(A) & = \sum_{i=1}^m\sum_{h=1}^ma_{ih} k_4(x^i,x^h) \\
 & = \sum_{i=1}^m\sum_{h=1}^m \bar\alpha_i\bar\alpha_h y^iy^h k_4(x^i,x^h) \geq 0,
\end{align*}
where the inequality follows from~\eqref{eq: nonegativity of F}. Thus, we showed that $F$ is \textcolor{revblue}{supermodular}.

Similarly, for $S\subseteq [n]$ and $e \in [n]\setminus S$, we can use~\eqref{eq:increase} and write
\begin{align*}
    \Delta_e(S) & = \sum_{i=1}^m\sum_{h=1}^m a_{ih} \left (\sum_{k=1}^{d}\binom{d}{k}\bigl(g^{ih}_e\bigr)^k\!\!{\sum_{j_1,\dots,j_{d-k}\in S}}\prod_{\ell=1}^{d-k}g^{ih}_{j_\ell} \right )\\
    & = \sum_{i=1}^m\sum_{h=1}^m \bar\alpha_i\bar\alpha_h y^iy^h k_5(x^i,x^h) \geq 0,
\end{align*}
where $k_5(\cdot,\cdot)$ is a kernel, and conclude that $F(S)$ is monotone non-decreasing. \end{proof}

In order to prove Proposition \ref{prop: H sub}, we first list a series of useful properties of submodular functions (see, e.g.,~\cite{krause2014submodular}). Recall that a submodular function $f$ defined on a set $E$ is said to be \textit{monotone nondecreasing} if for all $S\subseteq E, e \in E\setminus S$, it holds that \(f(S) \leq f(S\cup\{e\})\). Similarly, \(f\) is \textit{monotone nonincreasing} if for all $S\subseteq E, e \in E\setminus S$, we have \(f(S) \geq f(S\cup\{e\})\).

\begin{proposition}[\cite{krause2014submodular}]
Let \(f\) be a submodular (supermodular) function defined on a set \(E\). The following holds:
\begin{itemize}
    \item The \textit{opposite} function \(g\), defined as \(g(S) := -f(S)\) for all \(S \subseteq E\), is supermodular (submodular).
    \item The \textit{complement} function \(h\), defined as \(h(S) := f(E \setminus S)\) for all \(S \subseteq E\), is submodular (supermodular).
    \item The complement function of a monotone nondecreasing (nonincreasing) submodular function is a monotone nonincreasing (nondecreasing) submodular function.
\end{itemize}
\end{proposition}

\textbf{Proof of proposition \ref{prop: H sub}}

\begin{proof}

We recall that function $F$ is supermodular and nondecreasing. Take function $G = -F$, which is  submodular, nonincreasing. Problem \eqref{model: min supermod} is thus equivalent to:
\begin{equation}\label{model: max G}
\begin{array}{rl}
\max \limits_{ S} \; &  G(S)\\
\text{s.t.}\;  & S \in N\\
& |S| = \textcolor{revblue}{B}
\end{array}
\end{equation}

Consider function $\bar G$ such that $\bar G(S) = G([n] \setminus S) = -F([n] \setminus S) $, thus:
$$ \bar G(S) = \textcolor{revblue}{-}\sum_{i=1}^m\sum_{h=1}^m \bar\alpha_i\bar \alpha_h y^iy^h  \left (c + \sum_{j \in [n] \setminus S} \gamma x^i_j x^h_j \right )^d   .$$

Note that $\bar G$ is now submodular and  nondecreasing. Maximizing $G(S)$ over the constraint $|S| = B$, is equivalent to maximizing  $ \bar G(S)$ over the constraint $|S| = n-B$. Note that, because $ \bar G(S)$ is nondecreasing,  we can  maximize it over the constraint $| S| \leq n-B$ and obtain the same solution. Moreover, in order to have nonnegativity of the objective function, we can substract to  $\bar G$ its maximum value which is retained at the empty set $\emptyset$, thus obtaining the following optimization model which is equivalent to  \eqref{model: max G} and thus to \eqref{model: min supermod}:
\begin{equation}\label{model: proof reduction to max submod}
\begin{array}{rl}
\max \limits_{S} \; &\bar G( S) - \bar G(\emptyset) \\
\text{s.t.}\;  & S \in N \\
& |S| \leq n - B.
\end{array}
\end{equation}
Any \eqref{model: max G}-feasible solution $S'$ has a one-to-one correspondence to a  \eqref{model: proof reduction to max submod}-feasible solution solution $S''$, in the sense that $S' = [n] \setminus S''$ and $S'' = [n] \setminus S'$ . The objective function of \eqref{model: proof reduction to max submod} is submodular, nondecreasing and nonnegative, and  we have that
$$ \bar G(S) - \bar G(\emptyset)  = -  \sum_{i=1}^m\sum_{h=1}^m \bar\alpha_i\bar \alpha_h y^iy^h  \left (c + \sum_{j \in [n] \setminus S} \gamma x^i_j x^h_j \right )^d  + \underbrace{   \sum_{i=1}^m\sum_{h=1}^m \bar\alpha_i\bar \alpha_h y^iy^h  \left (c + \sum_{j \in [n]} \gamma x^i_j x^h_j \right )^d}_{\lambda}   = H(S)$$

This way, given a \eqref{model: proof reduction to max submod}-feasible solution solution $S''$, say set $S' = [n]\setminus S''$:

$$H(S'') = \bar G(S'') + \lambda = G(S') + \lambda = - F(S') + \lambda$$

and

$$ F(S') =  \lambda - H(S'')  $$

which concludes the proof.
\end{proof}

\subsection{Counterexamples to the super/submodularity of Gaussian Kernels}

We end the Appendix by showing that the function $F(S)$ is in general neither supermodular nor submodular when Gaussian kernels are used:
$$
F(S) = \sum_{i=1}^m \sum_{h=1}^m \bar\alpha_i \bar\alpha_h y^i y^h
\exp\left(-\gamma \sum_{j \in S}  (x^i_j - x^h_j)^2 \right).
$$

Consider the following dataset with three samples and parameters:
$$
\begin{aligned}
&x^1 = (0,\, 1,\,-0.5), \quad x^2 = (-1,\,-1,\,-1), \quad x^3 = (0,\,-1,\,0),\\
&y^1 = -1,\; y^2 =  y^3 = 1,\qquad
\bar\alpha_1 = 1,\; \bar\alpha_2 = \bar\alpha_3 = 0.5,\qquad \gamma = 1.
\end{aligned}
$$

For this setting, take first $A = \emptyset$, $B = \{1\}$, and $e = 2$. We compute:
$$
F(A) = 0,\quad F(A \cup \{e\}) \approx 1.9634,\quad
\Delta_e(A) \approx 1.9634,
$$
$$
F(B) \approx 0.3161,\quad F(B \cup \{e\}) \approx 1.6589,\quad
\Delta_e(B) \approx 1.3428.
$$
Since $\Delta_e(A) > \Delta_e(B)$, the supermodularity condition
$\Delta_e(A) \leq \Delta_e(B)$ is violated.

Similarly, take $A = \emptyset$, $B = \{1\}$, and $e = 3$. We compute:
$$
F(A) = 0,\quad F(A \cup \{e\}) \approx 0.1263,\quad
\Delta_e(A) \approx 0.1263,
$$
$$
F(B) \approx 0.3161,\quad F(B \cup \{e\}) \approx 0.5024,\quad
\Delta_e(B) \approx 0.1863.
$$
Since $\Delta_e(A) < \Delta_e(B)$, the submodularity condition
$\Delta_e(A) \geq \Delta_e(B)$ is violated.

\providecommand{\nref}{\#\text{SVM solves}}   %

\subsection{Configuration of the metaheuristic parameters}\label{app:param}
\color{revblue}

This section reports the study behind the parameter configurations adopted in
Section~\ref{sec: 5}. To make the conclusions robust, we repeat the analysis over \gnote{the nine datasets and} a grid of the four kernels (polynomial of
degree $d\in\{2,3,5\}$ with $c=1$, and the Gaussian kernel), five cardinality budgets
$B\in\{0.1\,n,0.2\,n,0.3\,n,0.4\,n,0.5\,n\}$, and the SVM hyper-parameters $C\in\{0.1,1,10,100\}$ and
$\gamma\in\{0.01,0.1,1\}$, for a total of $2160$ cases. Every configuration of parameters is treated as a
distinct algorithm, and all runs start from the same deterministic point, the optimal feature subset
obtained from the SVM dual computed on all features (the $\alpha$-start), so that we measure the effect of
the algorithm and not of a lucky initialization. We recall that \AOsub relies on the submodular
reformulation available only for polynomial kernels, so on the Gaussian kernel only \texttt{LS} and
\texttt{LS$^*$} apply.

Both heuristics parallelize naturally, and our implementation exploits this. At every iteration the dominant
cost is a large number of \emph{independent} convex SVM subproblems, one per candidate feature subset; we
dispatch these subproblems to a pool of worker threads, each running the \texttt{libsvm} solver, so that they
are solved in parallel and the wall-clock time decreases as the number of available cores grows. The
per-candidate objective-value computation $f(\alpha,\beta)$ (a dense linear-algebra operation) is kept
single-threaded, so that parallelism is exploited only across the independent subproblems and the cores are
not oversubscribed.

We measure quality as the gap to the best objective found by the two heuristic on each case. \gnote{Each of
these parameters can be increased indefinitely, and doing so buys a smaller and smaller quality
improvement at an ever growing cost. For every heuristic we therefore select, separately for each
parameter, the value beyond which a further increase only adds computational effort without improving the
solution quality: we treat the parameters one at a time because two values may reach the same objective
while differing severalfold in cost.} Alongside the running time, we report for each configuration the
total number of convex SVM subproblems it solves ($\nref$), a hardware-independent measure of
computational effort that does not depend on the degree of parallelism and enables a fair comparison
across machines.

\subsubsection*{Selecting the \texttt{LS$^*$} parameters}\label{sec:param:ls}

\lpnote{For \texttt{LS$^*$} we vary three parameters: the \emph{re-fit fraction}, the number $\#\mathrm{Samples}$ of candidates sampled from $\mathcal{N}^2$, and \texttt{OptWindow}, the parameter of the stopping criterion adopted in
Section~\ref{sec: 3}. \gnote{We consider the following ranges:}
\[
  \boxed{\;\texttt{LS$^*$}:\ \text{ re-fit}\in \{\text{single},5\%,10\%,50\%,100\%\},\ \#\mathrm{Samples}\in \{100,200,500,1000\},\ \texttt{OptWindow}\in \{5,10,15\}\;}
\]
}

\begin{figure}[!ht]\centering\color{revblue}
  \includegraphics[width=\linewidth]{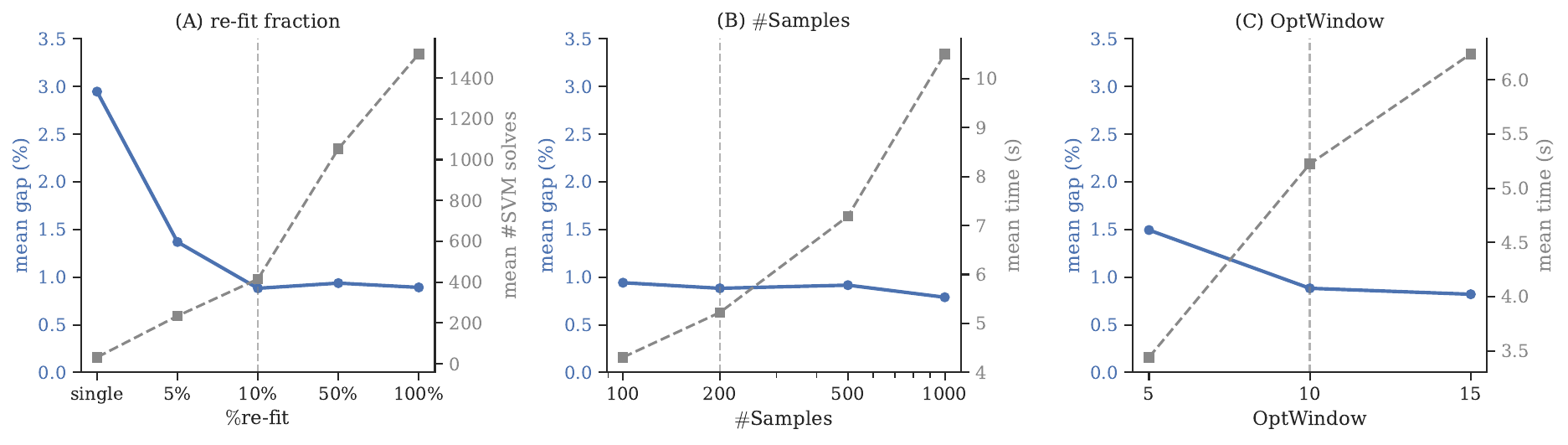}
  \caption{\color{revblue}{Parameter selection for \texttt{LS$^*$} (all $2160$ cases; blue = mean gap to the
  best local-search solution, grey = cost; the dashed vertical line marks the selected value). (A) re-fit
  fraction; (B) $\#$Samples; (C) \texttt{OptWindow}.}}
  \label{fig:param-ls}
\end{figure}

Figure~\ref{fig:param-ls} reports, for each of the three parameters, the mean gap to the best solution found together with the corresponding cost, over the $2160$ cases.

The \emph{re-fit fraction} is the most influential parameter: re-fitting the dual variables for
\emph{every} candidate move is accurate but expensive, so only a top fraction of the candidates (ranked by
a cheap screening score) is re-fitted exactly, and the rest are skipped. The quality curve has a sharp knee
at $10\%$ (Figure~\ref{fig:param-ls}(A)): the mean gap drops from $3.0\%$ (re-fitting only the single best
candidate) to $0.9\%$, and is flat from there on, while the cost keeps growing (mean \nref{} $415$ at
$10\%$ versus $1518$ at a full re-fit). \gnote{The $10\%$ setting therefore matches the quality of a full
re-fit at a fraction of the SVM solves.}

The number of sampled candidates, set by \texttt{$\#$Samples}, is by contrast far less sensitive
(Figure~\ref{fig:param-ls}(B)): going from $100$ to $1000$ sampled candidates moves the mean gap only from
$0.94\%$ to $0.79\%$, while the mean running time more than doubles (from $4.31$ to $10.51$ seconds). Note
that here the cost is \emph{screening time}, not SVM solves: the perturbations are only screened, so the
mean \nref{} stays at about $400$ throughout. We select $200$, at the beginning of the plateau.

Finally, \texttt{OptWindow} pays up to $10$ (Figure~\ref{fig:param-ls}(C)): raising it from $5$ to
$10$ lowers the mean gap from $1.49\%$ to $0.88\%$, whereas going from $10$ to $15$ \gnote{brings no
further quality gain} while the mean time keeps growing. On the Gaussian kernel, where only
\texttt{LS}/\texttt{LS$^*$} apply, these selected values remain the best trade-off, in particular the $10\%$
re-fit fraction.

\subsubsection*{Selecting the \AOsub parameters}\label{sec:param:ao}

\lpnote{For \AOsub we vary the pool size $\#\mathrm{SolPool}$ (how
many candidate subsets the submodular step keeps and hands to the exact SVM solve), and \texttt{OptWindow}, the parameter of the stopping criterion adopted in Section~\ref{sec: 4}; the basic scheme \AObase{} is
parameter-free. \gnote{We consider the following ranges:}
\[
  \boxed{\;\AOsub:\ \#\mathrm{SolPool}\in \{10,50,100,200\},\ \texttt{OptWindow}\in \{5,10,15\}\;}
\]}

\begin{figure}[!ht]\centering\color{revblue}
  \includegraphics[width=0.8\linewidth]{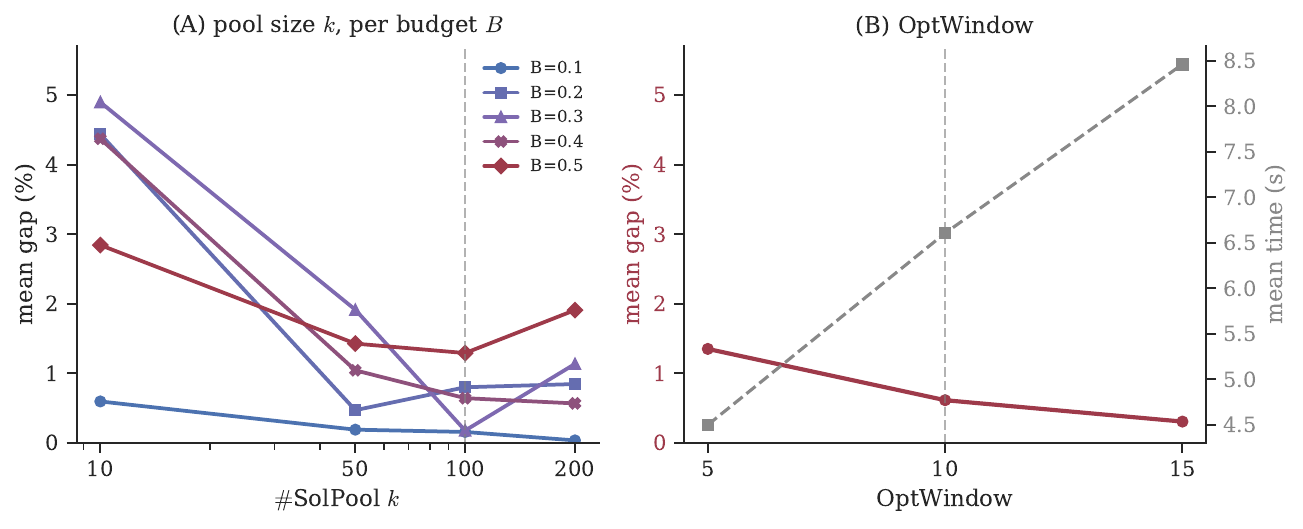}
  \caption{\color{revblue}Parameter selection for \AOsub (polynomial cases; the dashed vertical line marks
  the selected value). (A) pool size, one curve per budget $B$; (B) \texttt{OptWindow}.}
  \label{fig:param-ao}
\end{figure}

Figure~\ref{fig:param-ao} reports the same analysis for the two parameters of \AOsub, on the $1620$
polynomial cases, the gap being measured here against the best solution found.

The pool size \texttt{$\#\mathrm{SolPool}$}, unlike the far less sensitive parameters of \texttt{LS$^*$}, is
genuinely budget-dependent (Figure~\ref{fig:param-ao}(A), one curve per budget). A small pool ($k{=}50$)
is enough at small budgets but leaves a visible gap at the harder mid-range budget $B{=}0.3\,n$ (mean gap
$2.1\%$), which $k{=}100$ closes ($0.2\%$) at essentially the same cost: a larger budget means a larger
combinatorial space for the submodular $\beta$-step, so more candidate solutions must be kept. Enlarging
the pool further does \emph{not} pay: at $k{=}200$ several curves bend upwards (mean gap from $0.69\%$
back to $1.01\%$), because with a larger pool the re-computation of $\alpha$ more often moves the search
to suboptimal $\beta$ candidates (cf.\ the remark at the end of Section~\ref{sec: 4}). \gnote{The pool
$k{=}100$ gives the best mean quality at a moderate cost, and is our choice.}

As for \texttt{OptWindow}, the alternating scheme benefits from a longer window slightly more than the local search
(Figure~\ref{fig:param-ao}(B)): raising it from $5$ to $10$ lowers the mean gap
from $1.35\%$ to $0.61\%$, and from $10$ to $15$ to $0.31\%$. \gnote{The residual gain at $15$
is, however, paid with a further increase in running time}, so we keep $10$ as
the best quality--time compromise.

\subsection{Comparison of the heuristics}\label{app:compare}

\gnote{Together with the two improved heuristics used throughout the paper, \texttt{LS$^*$} and \AOsub,
we compare here also their basic versions. By \texttt{LS} we denote the
plain local search of Section~\ref{sec: 3}, which stops as soon as no single swap improves the objective,
without the sampling of $\mathcal{N}^2$ of \texttt{LS$^*$}; by \AObase{} we denote the basic alternating scheme of
Section~\ref{sec: 4}, which stops at the first iteration that does not improve the objective, without the
tabu list and the pool of candidate solutions of \AOsub. Comparing each enhanced version with its basic
counterpart isolates the contribution of the enhancements. All the comparisons reported here are run on
the same grid of Section~\ref{app:param}, restricted to its $1620$ polynomial cases, $324$ for each
budget: \AObase{} and \AOsub rely on the submodular reformulation and are therefore not defined on the
Gaussian kernel. Throughout this section the gap of a method on a given case is measured against the best
solution found on that case by any of the four algorithms, and is then averaged over the cases.}

\begin{figure}[!ht]\centering\color{revblue}
  \includegraphics[width=0.85\linewidth]{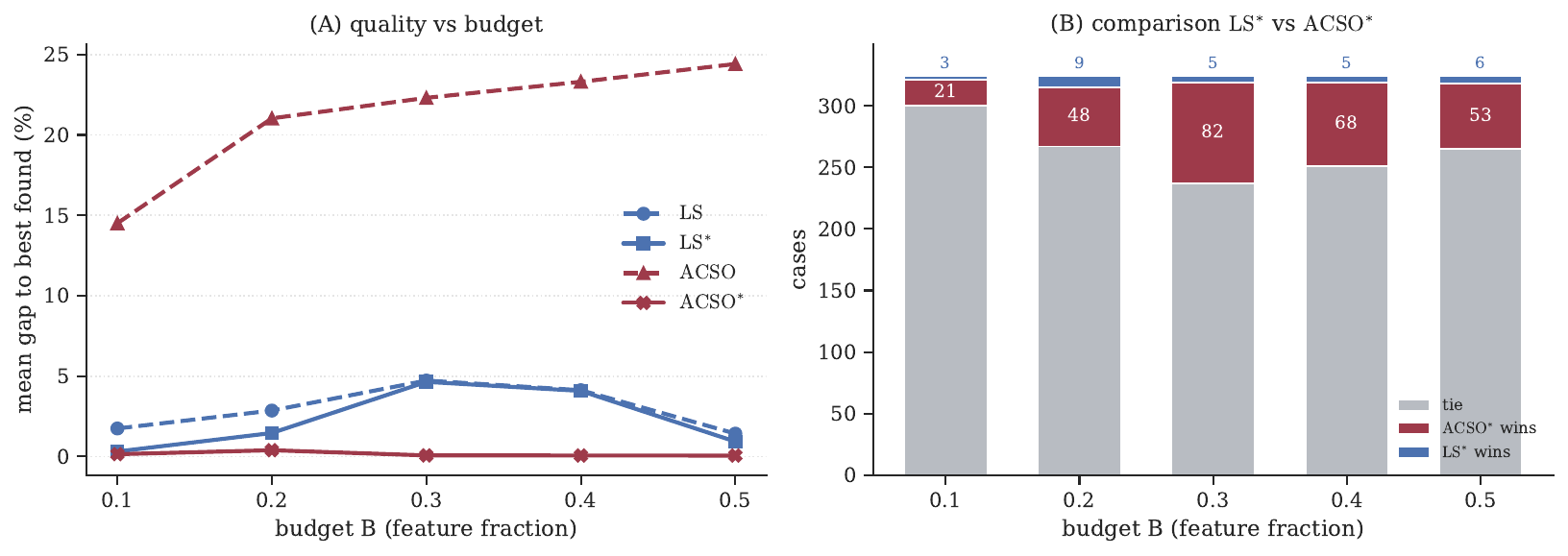}
  \caption{\textcolor{revblue}{Effect of the budget $B$ (polynomial kernels, averaged over datasets, $C$, $d$
  and $\gamma$). (A) mean gap to the best found solution for the four heuristics. (B) direct comparison of
  the selected \texttt{LS$^*$} and \AOsub, per budget ($324$ cases): ties in grey, wins (strictly better
  objective) in colour.}}
  \label{fig:gapB}
\end{figure}

\begin{figure}[!ht]\centering\color{revblue}
  \includegraphics[width=0.8\linewidth]{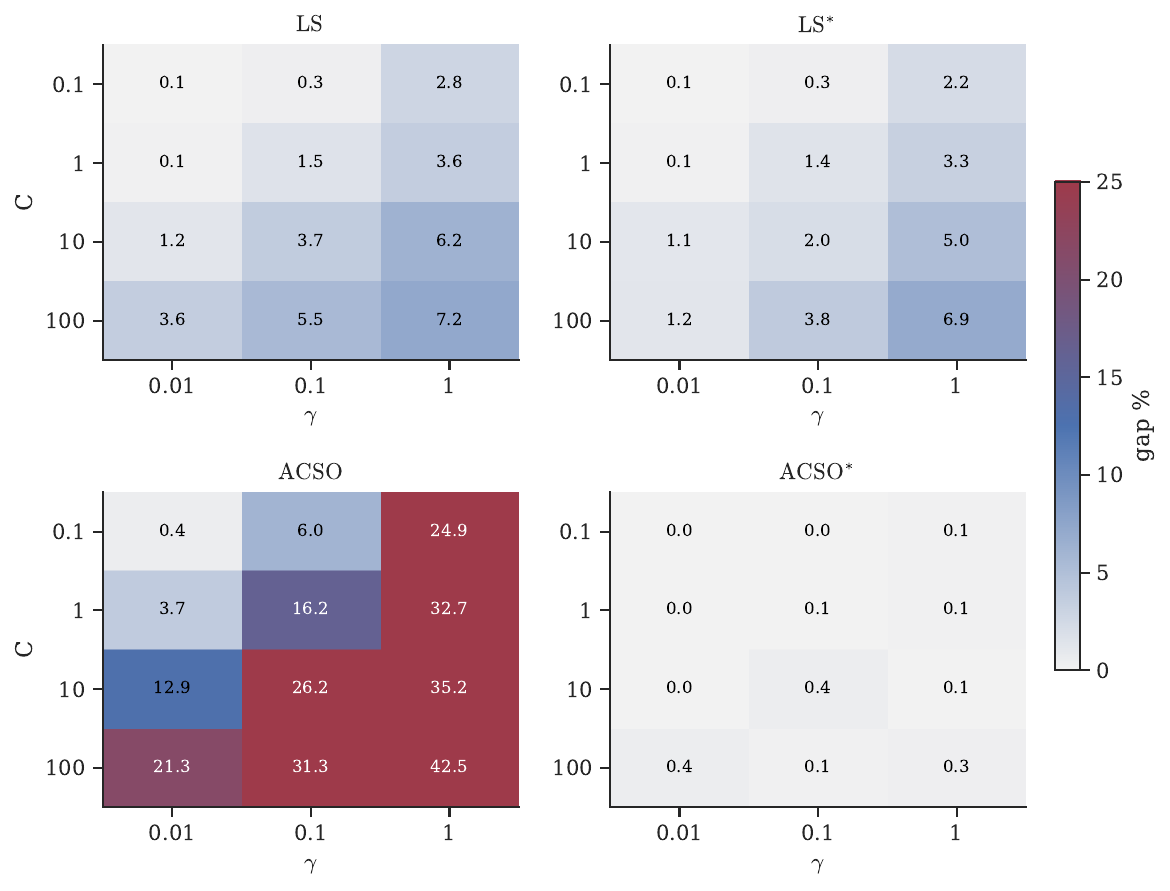}
  \caption{\textcolor{revblue}{Mean gap to the best found solution over the $C$--$\gamma$ grid (polynomial
  cases, averaged over datasets and budgets), one panel per method.}}
  \label{fig:cgamma}
\end{figure}

\gnote{The comparison of the four heuristics is reported from three complementary angles.
Figure~\ref{fig:gapB}(A) shows how their mean gap evolves as the budget $B$ grows, averaged over
datasets, polynomial kernels and the $C$--$\gamma$ grid; panel~(B) of the same figure, discussed further
below, compares the two improved versions head to head. Table~\ref{tab:ladder} breaks the same quantity
down per dataset, averaging over the polynomial kernels, the $C$--$\gamma$ grid and the five budgets.
Figure~\ref{fig:cgamma} resolves it cell by cell over the $C$--$\gamma$ grid, averaging over datasets and
budgets, with one panel per heuristic.}

\begin{table}[!ht]
\centering\color{revblue}\renewcommand\arraystretch{1.2}
\begin{tabular}{lc|cc|cc}
\toprule
dataset & $n$ & \texttt{LS} & \texttt{LS$^*$} & \AObase & \AOsub \\
\midrule
Wholesale  & 7  & $\mathbf{0.00}$ & $\mathbf{0.00}$ & $13.11$ & $\mathbf{0.00}$ \\
Diabetes   & 8  & $0.06$ & $0.05$ & $4.78$  & $\mathbf{0.00}$ \\
BCW        & 9  & $0.50$ & $0.20$ & $30.60$ & $\mathbf{0.00}$ \\
Cleveland  & 13 & $0.10$ & $0.04$ & $17.03$ & $\mathbf{0.00}$ \\
Parkinsons & 22 & $0.66$ & $0.85$ & $9.42$  & $\mathbf{0.03}$ \\
German     & 24 & $0.07$ & $\mathbf{0.06}$ & $5.39$  & $\mathbf{0.06}$ \\
BCD        & 30 & $12.46$ & $8.57$ & $41.78$ & $\mathbf{0.57}$ \\
Ionosphere & 33 & $4.22$ & $3.15$ & $32.44$ & $\mathbf{0.11}$ \\
Sonar      & 60 & $8.73$ & $7.65$ & $35.45$ & $\mathbf{0.56}$ \\
\bottomrule
\end{tabular}
\caption{\textcolor{revblue}{Mean gap to the best found solution per dataset (averaged over polynomial degrees, $C$, $\gamma$,
and the five budgets). \texttt{LS} = plain local search; \texttt{LS$^*$} = $10\%$ re-fit; \AObase = basic
alternating scheme (no tabu, single candidate); \AOsub = $k{=}100$. Bold marks the best value in each row.}}
\label{tab:ladder}
\end{table}

\gnote{However the results are aggregated, the picture is the same. The two improved heuristics stay
close to the best found solution throughout, and \AOsub is the more stable of them: its mean gap never
exceeds $0.4\%$ in any cell of the $C$--$\gamma$ grid, and stays within $0.6\%$ on every dataset. The two
basic versions, on the contrary, are never competitive. \AObase{} degrades from $15\%$ to $24\%$ as the
budget grows (Figure~\ref{fig:gapB}(A)), is far from the best on \emph{every} dataset, from $5\%$ up to
$42\%$ (Table~\ref{tab:ladder}), and ranges from $0.4\%$ in the mildest cell of the grid to $42.5\%$ in
the harshest (Figure~\ref{fig:cgamma}).} The decomposition alone is therefore not enough, and it is the
tabu-and-pool enhancement that gives \AOsub its quality, just as the sampling of $\mathcal{N}^2$ and the
tabu list give \texttt{LS$^*$} its own. The plain local search \texttt{LS} finds the best solution on the
small datasets but leaves a large gap on the high-dimensional ones, BCD ($12.5\%$), Sonar ($8.7\%$),
Ionosphere ($4.2\%$), and the richer re-fit of \texttt{LS$^*$} only partly closes it. \gnote{The
difficulty is driven by the number of features: the three datasets with $n\ge 30$ are the only ones on
which the local searches leave a sizeable gap, whereas below that threshold they stay under $0.7\%$. All
four keep their relative standing over the whole $C$--$\gamma$ grid, the two basic versions degrading
fastest where $C$ and $\gamma$ are largest, the regime in which the kernel matrix is most ill-conditioned
and the problem hardest.}

\gnote{We now turn to the two improved heuristics alone. Figure~\ref{fig:gapB}(B) counts, for each
budget, on how many of the $324$ polynomial cases they reach the same objective and on how many each of
them is strictly better than the other. Table~\ref{tab:quality-cost} completes the comparison with their
cost, reporting per budget the mean gap together with the mean number of SVM subproblems solved and the
mean running time.}

\begin{table}[!ht]
\centering\color{revblue}\renewcommand\arraystretch{1.2}
\begin{tabular}{c|cc|cc|cc}
\toprule
& \multicolumn{2}{c|}{mean gap (\%)} & \multicolumn{2}{c|}{mean $\nref$} & \multicolumn{2}{c}{mean time (s)} \\
$B$ & \texttt{LS$^*$} & \AOsub & \texttt{LS$^*$} & \AOsub & \texttt{LS$^*$} & \AOsub \\
\midrule
$0.1\,n$ & $0.31$ & $\mathbf{0.15}$ & $\mathbf{178}$ & $362$ & $\mathbf{2.51}$ & $3.33$ \\
$0.2\,n$ & $1.46$ & $\mathbf{0.39}$ & $\mathbf{360}$ & $655$ & $\mathbf{4.21}$ & $5.74$ \\
$0.3\,n$ & $4.64$ & $\mathbf{0.07}$ & $\mathbf{463}$ & $783$ & $\mathbf{5.84}$ & $8.03$ \\
$0.4\,n$ & $4.08$ & $\mathbf{0.07}$ & $\mathbf{532}$ & $771$ & $\mathbf{7.13}$ & $8.23$ \\
$0.5\,n$ & $0.93$ & $\mathbf{0.06}$ & $\mathbf{570}$ & $696$ & $\mathbf{7.10}$ & $7.70$ \\
\bottomrule
\end{tabular}
\caption{\textcolor{revblue}{Quality and cost of the two selected configurations, per budget, over the $324$ polynomial cases
of each budget: mean gap to the best found solution, mean number of SVM subproblems solved, and mean
running time. Bold marks the better value of each pair.}}
\label{tab:quality-cost}
\end{table}

The two are \emph{equivalent on most cases}: they reach the same objective in $237$ to $300$ of the $324$
polynomial cases per budget. Where they differ, \AOsub almost always wins, and its wins concentrate at
the mid-range budgets (up to $82$ against $5$ at $B{=}0.3\,n$, the combinatorially hardest point).
\AOsub{} also attains the smaller gap at every budget, whereas \texttt{LS$^*$} has the lower cost on both
measures: it solves fewer SVM subproblems at every budget, since most of its candidate moves are only
screened rather than re-solved exactly, and it also has the smaller mean running time at every budget.
The cost of \AOsub{} is also the more dispersed of the two across cases. The trade-off between the two
heuristics is therefore between solution quality and running cost, with the further difference that
\texttt{LS$^*$} applies to any kernel, \AOsub{} being restricted to the polynomial case.

\normalcolor

\subsection{\textcolor{revblue}{Optimality certification procedure}}\label{app:cert}
{\color{revblue}
For a fixed feature vector $\beta$ the problem \eqref{prob:DFS-SVM} reduces to
a \emph{convex} SVM dual. We can therefore enumerate \emph{all} $ \binom{n}{B}$ cardinality-$B$ feature subsets, solve the
corresponding convex dual exactly (using e.g. LIBSVM, through \texttt{scikit-learn}), and take the smallest objective value. This value is the certified global optimum and the minimizing subset a globally optimal selection. The
procedure is exact and independent of the kernel degree.
Thus, the hardness of the certification has two sources: the combinatorial number $\binom{n}{B}$ of
cardinality-$B$ subsets to enumerate, and the cost of each individual convex SVM re-fit, which grows with
the number of samples $m$. To keep the enumeration tractable, we restrict the budget to
$B\in[1,\lfloor n/2\rfloor]$ and, within this range, stop at the largest value $B_{\max}$ for which both
the number of subsets and the total re-fitting cost remain affordable; in practice this means
$\binom{n}{B_{\max}}\le 10^{8}$, \gnote{and a smaller $B_{\max}$ on the datasets with a large number of
samples $m$, where each individual convex subproblem is more expensive to solve}. We cap $B$ at $\lfloor n/2\rfloor$
because $\binom{n}{B}=\binom{n}{n-B}$: budgets beyond $\lfloor n/2\rfloor$ mirror the combinatorial
difficulty of smaller ones, so enumerating them would merely duplicate combinatorially easy instances and
optimistically bias the aggregate results. Restricting to $[1,\lfloor n/2\rfloor]$ instead yields a clean,
monotonically increasing difficulty scale, from the easiest budget $B=1$ to the hardest
$B=\lfloor n/2\rfloor$. The resulting $B_{\max}$ is reported per dataset in
Table~\ref{tab:cert-coverage}: \gnote{it equals $\lfloor n/2\rfloor$ on the datasets with few features,
and it is smaller on those with many features or many samples} (down to $B_{\max}=6$ for Sonar, where
$\binom{n}{B}$ would otherwise reach about $10^{17}$), beyond which exact certification becomes
intractable even on a high-performance computer.
The underlined values in Table~\ref{tab:cert-coverage} are those for which $B_{\max}<\lfloor n/2\rfloor$:
for Breast Cancer Diagnostic, Ionosphere and Sonar the binding limit is the number of combinations
mentioned above, whereas for German the certification stops well before that limit is approached,
\gnote{because with $m=1000$ samples each convex subproblem is far more expensive to solve}.

Table~\ref{tab:cert-coverage} summarizes the instances solved to optimality.

\begin{table}[!ht]
\centering\color{revblue}\renewcommand\arraystretch{1.2}\footnotesize
\begin{tabular}{lrrcrr}
\toprule
dataset & $n$ & $m$ & $B_{\max}$ & $\binom{n}{B_{\max}}$ & \# QP solved \\
\midrule
Wholesale  &  7 & 440 & $3$  & $35$              & $2.5\times10^{2}$ \\
Diabetes   &  8 & 768 & $4$  & $70$              & $6.5\times10^{2}$ \\
BCW &  9 & 699 & $4$  & $126$             & $1.0\times10^{3}$ \\
Cleveland  & 13 & 297 & $6$  & $1.7\times10^{3}$ & $1.6\times10^{4}$ \\
Parkinsons & 22 & 195 & $11$ & $7.1\times10^{5}$ & $9.8\times10^{6}$ \\
German     & 24 & 1000 & $\underline{9}$  & $1.3\times10^{6}$ & $1.0\times10^{7}$ \\
BCD   & 30 & 569 & $\underline{12}$ & $8.6\times10^{7}$ & $6.9\times10^{8}$ \\
Ionosphere & 33 & 351 & $\underline{10}$ & $9.3\times10^{7}$ & $5.1\times10^{8}$ \\
Sonar      & 60 & 208 & $\underline{6}$  & $5.0\times10^{7}$ & $2.2\times10^{8}$ \\
\bottomrule
\end{tabular}
\caption{%
\lpnote{
For each dataset: number of features $n$, number of samples
$m$, $B_{\max}$ for which  optima are certified, size of the largest enumeration $\binom{n}{B_{\max}}$, and total
number of convex quadratic programs solved (summed over the polynomial kernels with degrees $d=2,3,5$ and the Gaussian
kernel, where certified). Underlined values are those for which $B_{\max}<\lfloor n/2\rfloor$.}
}
\label{tab:cert-coverage}
\end{table}

}

\subsection{Comparison with \texttt{Baron} on reduced datasets}\label{sec:baron}
Baron is one of the best exact commercial solvers for general MINLPs, nevertheless for datasets with more than 200 samples the solver is not able to find a solution in a two hours timelimit. We decided to compare the results on different instances created by reducing the size of four datasets, Cleveland, Breast Cancer Diagnostic, Ionosphere and Sonar. The Cleveland dataset is representative of the behaviour of the algorithms on the datasets with fewer features. Breast Cancer Diagnostic, Ionosphere and Sonar datasets are instead much harder to solve having more than 30 features. \textcolor{revblue}{On the same reduced instances we additionally run the general-purpose nonlinear solvers \texttt{Knitro} and \texttt{Ipopt} (described in Section~\ref{sec: 5}), so that all competing methods are evaluated on a common set of instances. For a fair comparison, \texttt{Knitro} and \texttt{Ipopt} are run from the same five multistart points used by \AOsub and \texttt{LS$^*$}, and we report, as for our methods, the best objective value found together with the cumulative running time over the five starts.}

{In these results, we refer to \texttt{Baron} as the complete optimization process of the solver, including both its initial heuristic and local search procedures, as well as the global branch-and-bound phase. In contrast, \texttt{Baron-H} denotes the solution obtained by terminating \texttt{Baron} after this initial phase, before the global optimization begins. We report \texttt{Baron-H} separately because, for hard instances like those considered, \texttt{Baron} often fails to prove optimality, while its initial phase alone typically provides the best available solution.}

\textcolor{revblue}{For the large majority of these instances \texttt{Baron} cannot close the optimality gap within the time limit, so the MIP gaps it returns are uninformative and do not certify the quality of the reported solutions: the only instance solved to certified optimality is Cleveland with $d=5$ and $m=10$.}

\begin{table}[!ht]
\centering
\renewcommand\arraystretch{1.2}\footnotesize
\resizebox{\textwidth}{!}{%
\begin{tabular}{ll|cc|cc|cc|cc|cc|cc}
\toprule
 & & \multicolumn{2}{c|}{\texttt{Baron}} & \multicolumn{2}{c|}{\texttt{Baron-H}} & \multicolumn{2}{c|}{\textcolor{revblue}{\texttt{Knitro}}} & \multicolumn{2}{c|}{\textcolor{revblue}{\texttt{Ipopt}}} & \multicolumn{2}{c|}{\AOsub} & \multicolumn{2}{c}{\texttt{LS$^*$}}\\
 & & ObjFun & MipGap(\%)  & ObjFun & Time & \textcolor{revblue}{ObjFun} & \textcolor{revblue}{Time} & \textcolor{revblue}{ObjFun} & \textcolor{revblue}{Time} & ObjFun & \textcolor{revblue}{Time} & ObjFun & \textcolor{revblue}{Time} \\
\midrule
\multirow{5}{*}{\texttt{d=2}} & \texttt{m=10} & \textbf{-41.1} & 5.0e0 & \textbf{-41.1} & \textbf{0.1} & \textcolor{revblue}{\textbf{-41.1}} & \textcolor{revblue}{0.3} & \textcolor{revblue}{-38.8} & \textcolor{revblue}{0.2} & \textbf{-41.1} & \textcolor{revblue}{0.3} & \textbf{-41.1} & \textcolor{revblue}{0.8} \\
 & \texttt{m=50} & -227.8 & 1.72e4  & -224.2 & 1.8 & \textcolor{revblue}{\textbf{-228.1}} & \textcolor{revblue}{2.5} & \textcolor{revblue}{\textbf{-228.1}} & \textcolor{revblue}{2.2} & \textbf{-228.1} & \textcolor{revblue}{0.8} & \textbf{-228.1} & \textcolor{revblue}{\textbf{0.7}} \\
 & \texttt{m=100} & \textbf{-625.0} & 3.36e4 & \textbf{-625.0} & 19.4 & \textcolor{revblue}{-565.3} & \textcolor{revblue}{9.6} & \textcolor{revblue}{-611.3} & \textcolor{revblue}{10.2} & \textbf{-625.0} & \textcolor{revblue}{\textbf{0.9}} & \textbf{-625.0} & \textcolor{revblue}{1.0} \\
 & \texttt{m=150} & \textbf{-930.9} & 4.98e4  & \textbf{-930.9} & 91.4 & \textcolor{revblue}{-853.2} & \textcolor{revblue}{23.4} & \textcolor{revblue}{\textbf{-930.9}} & \textcolor{revblue}{26.8} & \textbf{-930.9} & \textcolor{revblue}{\textbf{1.2}} & \textbf{-930.9} & \textcolor{revblue}{1.5} \\
 & \texttt{m=200} & -1267.6 & 1.59e6  & -1267.6 & 286.9 & \textcolor{revblue}{\textbf{-1286.5}} & \textcolor{revblue}{42.3} & \textcolor{revblue}{-1267.6} & \textcolor{revblue}{55.3} & \textbf{-1286.5} & \textcolor{revblue}{\textbf{1.5}} & \textbf{-1286.5} & \textcolor{revblue}{1.7} \\
\midrule
\multirow{5}{*}{\texttt{d=3}} & \texttt{m=10} & \textbf{-26.1} & 5.0e0  & -20.9 & \textbf{0.1} & \textcolor{revblue}{-24.7} & \textcolor{revblue}{0.2} & \textcolor{revblue}{\textbf{-26.1}} & \textcolor{revblue}{0.2} & \textbf{-26.1} & \textcolor{revblue}{0.5} & \textbf{-26.1} & \textcolor{revblue}{0.7} \\
 & \texttt{m=50} & \textbf{-191.7} & 2.70e4  & \textbf{-191.7} & 1.8 & \textcolor{revblue}{-165.0} & \textcolor{revblue}{3.0} & \textcolor{revblue}{-165.0} & \textcolor{revblue}{2.8} & \textbf{-191.7} & \textcolor{revblue}{0.9} & \textbf{-191.7} & \textcolor{revblue}{\textbf{0.9}} \\
 & \texttt{m=100} & \textbf{-468.1} & 7.62e4  & \textbf{-468.1} & 18.8 & \textcolor{revblue}{-454.6} & \textcolor{revblue}{9.9} & \textcolor{revblue}{-454.3} & \textcolor{revblue}{10.3} & \textbf{-468.1} & \textcolor{revblue}{\textbf{1.0}} & \textbf{-468.1} & \textcolor{revblue}{1.1} \\
 & \texttt{m=150} & \textbf{-747.4} & 1.06e5 & \textbf{-747.4} & 91.9 & \textcolor{revblue}{-662.6} & \textcolor{revblue}{23.3} & \textcolor{revblue}{\textbf{-747.4}} & \textcolor{revblue}{26.9} & \textbf{-747.4} & \textcolor{revblue}{\textbf{1.3}} & \textbf{-747.4} & \textcolor{revblue}{1.5} \\
 & \texttt{m=200} & -1038.0 & 2.90e6  & -1038.0 & 295.4 & \textcolor{revblue}{-938.2} & \textcolor{revblue}{43.9} & \textcolor{revblue}{-1038.0} & \textcolor{revblue}{55.2} & \textbf{-1097.8} & \textcolor{revblue}{\textbf{1.4}} & \textbf{-1097.8} & \textcolor{revblue}{2.0} \\
\midrule
\multirow{5}{*}{\texttt{d=5}} & \texttt{m=10} & \textbf{-20.6} & 0.0 (933.2\emph{s}) & \textbf{-20.6} & \textbf{0.1} & \textcolor{revblue}{\textbf{-20.6}} & \textcolor{revblue}{0.2} & \textcolor{revblue}{\textbf{-20.6}} & \textcolor{revblue}{0.2} & \textbf{-20.6} & \textcolor{revblue}{0.7} & \textbf{-20.6} & \textcolor{revblue}{0.7} \\
 & \texttt{m=50} & -118.2 & 1.19e5  & -96.7 & 2.1 & \textcolor{revblue}{-119.2} & \textcolor{revblue}{2.9} & \textcolor{revblue}{-92.4} & \textcolor{revblue}{2.9} & \textbf{-149.3} & \textcolor{revblue}{0.9} & \textbf{-149.3} & \textcolor{revblue}{\textbf{0.8}} \\
 & \texttt{m=100} & -298.9 & 3.09e5  & -298.9 & 20.2 & \textcolor{revblue}{-245.7} & \textcolor{revblue}{9.7} & \textcolor{revblue}{-311.9} & \textcolor{revblue}{10.2} & \textbf{-347.1} & \textcolor{revblue}{\textbf{1.0}} & -327.2 & \textcolor{revblue}{1.7} \\
 & \texttt{m=150} & -458.0 & 4.74e6  & -458.0 & 90.5 & \textcolor{revblue}{-422.5} & \textcolor{revblue}{23.4} & \textcolor{revblue}{-486.9} & \textcolor{revblue}{29.0} & \textbf{-555.8} & \textcolor{revblue}{\textbf{1.2}} & -532.5 & \textcolor{revblue}{2.5} \\
 & \texttt{m=200} & -728.3 & 5.35e6  & -728.3 & 285.0 & \textcolor{revblue}{-634.4} & \textcolor{revblue}{42.7} & \textcolor{revblue}{-640.4} & \textcolor{revblue}{52.0} & \textbf{-824.5} & \textcolor{revblue}{\textbf{1.5}} & \textbf{-824.5} & \textcolor{revblue}{2.8} \\
\bottomrule
\end{tabular}}
\caption{Results for \texttt{Baron} (time limit 3600 sec), \texttt{Baron-H}, \textcolor{revblue}{\texttt{Knitro}, \texttt{Ipopt},} \AOsub and \texttt{LS*} on polynomial kernels - Cleveland dataset}
\label{table: baron results clev}
\end{table}

\begin{table}[!ht]
\centering
\renewcommand\arraystretch{1.2}\footnotesize
\resizebox{\textwidth}{!}{%
\begin{tabular}{ll|cc|cc|cc|cc|cc|cc}
\toprule
 & & \multicolumn{2}{c|}{\texttt{Baron}} & \multicolumn{2}{c|}{\texttt{Baron-H}} & \multicolumn{2}{c|}{\textcolor{revblue}{\texttt{Knitro}}} & \multicolumn{2}{c|}{\textcolor{revblue}{\texttt{Ipopt}}} & \multicolumn{2}{c|}{\AOsub} & \multicolumn{2}{c}{\texttt{LS$^*$}}\\
 & & ObjFun & MipGap(\%)  & ObjFun & Time & \textcolor{revblue}{ObjFun} & \textcolor{revblue}{Time} & \textcolor{revblue}{ObjFun} & \textcolor{revblue}{Time} & ObjFun & \textcolor{revblue}{Time} & ObjFun & \textcolor{revblue}{Time} \\
\midrule
\multirow{5}{*}{\texttt{d=2}} & \texttt{m=10} & \textbf{-2.13} & 1.80e4 & -2.08 & \textbf{0.1} & \textcolor{revblue}{-2.08} & \textcolor{revblue}{0.3} & \textcolor{revblue}{-2.08} & \textcolor{revblue}{0.2} & \textbf{-2.13} & \textcolor{revblue}{0.4} & \textbf{-2.13} & \textcolor{revblue}{1.7} \\
 & \texttt{m=50} & \textbf{-51.32} & 1.34e5 & \textbf{-51.32} & 4.0 & \textcolor{revblue}{-48.79} & \textcolor{revblue}{4.8} & \textcolor{revblue}{\textbf{-51.32}} & \textcolor{revblue}{4.8} & \textbf{-51.32} & \textcolor{revblue}{2.5} & \textbf{-51.32} & \textcolor{revblue}{\textbf{1.9}} \\
 & \texttt{m=100} & -89.23 & 3.50e5 & -89.23 & 46.6 & \textcolor{revblue}{-87.44} & \textcolor{revblue}{19.7} & \textcolor{revblue}{-89.43} & \textcolor{revblue}{17.3} & \textbf{-94.99} & \textcolor{revblue}{2.6} & \textbf{-94.99} & \textcolor{revblue}{\textbf{2.1}} \\
 & \texttt{m=150} & -122.56 & 6.20e5 & -122.56 & 213.8 & \textcolor{revblue}{-120.10} & \textcolor{revblue}{44.7} & \textcolor{revblue}{-125.55} & \textcolor{revblue}{46.9} & \textbf{-142.22} & \textcolor{revblue}{2.9} & \textbf{-142.22} & \textcolor{revblue}{\textbf{2.3}} \\
 & \texttt{m=200} & -160.67 & 7.69e6 & -160.67 & 698.9 & \textcolor{revblue}{-116.94} & \textcolor{revblue}{82.2} & \textcolor{revblue}{-178.05} & \textcolor{revblue}{102.3} & \textbf{-183.58} & \textcolor{revblue}{6.0} & \textbf{-183.58} & \textcolor{revblue}{\textbf{3.1}} \\
\midrule
\multirow{5}{*}{\texttt{d=3}} & \texttt{m=10} & \textbf{-1.06} & 2.16e2 & \textbf{-1.06} & \textbf{0.1} & \textcolor{revblue}{\textbf{-1.06}} & \textcolor{revblue}{0.3} & \textcolor{revblue}{\textbf{-1.06}} & \textcolor{revblue}{1.0} & \textbf{-1.06} & \textcolor{revblue}{0.7} & \textbf{-1.06} & \textcolor{revblue}{1.8} \\
 & \texttt{m=50} & -33.21 & 2.71e5 & -18.74 & 4.2 & \textcolor{revblue}{-31.47} & \textcolor{revblue}{5.4} & \textcolor{revblue}{\textbf{-34.08}} & \textcolor{revblue}{3.8} & \textbf{-34.08} & \textcolor{revblue}{2.5} & \textbf{-34.08} & \textcolor{revblue}{\textbf{2.0}} \\
 & \texttt{m=100} & -59.33 & 9.86e5 & -59.33 & 45.0 & \textcolor{revblue}{-53.21} & \textcolor{revblue}{18.9} & \textcolor{revblue}{-66.57} & \textcolor{revblue}{18.3} & \textbf{-71.00} & \textcolor{revblue}{\textbf{1.9}} & \textbf{-71.00} & \textcolor{revblue}{2.3} \\
 & \texttt{m=150} & -67.77 & 2.64e7 & -67.77 & 209.4 & \textcolor{revblue}{-85.38} & \textcolor{revblue}{44.0} & \textcolor{revblue}{-63.88} & \textcolor{revblue}{48.3} & \textbf{-112.78} & \textcolor{revblue}{5.1} & \textbf{-112.78} & \textcolor{revblue}{\textbf{2.2}} \\
 & \texttt{m=200} & -112.71 & 2.66e7 & -112.70 & 945.5 & \textcolor{revblue}{-98.79} & \textcolor{revblue}{84.1} & \textcolor{revblue}{-108.24} & \textcolor{revblue}{103.7} & -112.70 & \textcolor{revblue}{4.3} & \textbf{-117.63} & \textcolor{revblue}{\textbf{2.6}} \\
\midrule
\multirow{5}{*}{\texttt{d=5}} & \texttt{m=10} & \textbf{-0.36} & 1.28e4 & \textbf{-0.36} & \textbf{0.1} & \textcolor{revblue}{\textbf{-0.36}} & \textcolor{revblue}{0.4} & \textcolor{revblue}{\textbf{-0.36}} & \textcolor{revblue}{0.2} & \textbf{-0.36} & \textcolor{revblue}{0.7} & \textbf{-0.36} & \textcolor{revblue}{1.7} \\
 & \texttt{m=50} & \textbf{-17.07} & 2.08e6 & \textbf{-17.07} & 4.2 & \textcolor{revblue}{\textbf{-17.07}} & \textcolor{revblue}{5.3} & \textcolor{revblue}{\textbf{-17.07}} & \textcolor{revblue}{3.8} & \textbf{-17.07} & \textcolor{revblue}{2.7} & \textbf{-17.07} & \textcolor{revblue}{\textbf{2.2}} \\
 & \texttt{m=100} & \textbf{-29.34} & 3.32e8 & \textbf{-29.34} & 47.1 & \textcolor{revblue}{-22.44} & \textcolor{revblue}{18.9} & \textcolor{revblue}{-23.32} & \textcolor{revblue}{18.7} & -23.58 & \textcolor{revblue}{\textbf{2.1}} & \textbf{-29.34} & \textcolor{revblue}{\textbf{2.1}} \\
 & \texttt{m=150} & -32.66 & 1.24e9 & -32.66 & 217.0 & \textcolor{revblue}{-30.96} & \textcolor{revblue}{44.2} & \textcolor{revblue}{-32.53} & \textcolor{revblue}{53.0} & \textbf{-48.62} & \textcolor{revblue}{3.3} & \textbf{-48.62} & \textcolor{revblue}{\textbf{2.2}} \\
 & \texttt{m=200} & -60.65 & 1.11e9 & -60.65 & 742.1 & \textcolor{revblue}{-45.98} & \textcolor{revblue}{82.1} & \textcolor{revblue}{-60.17} & \textcolor{revblue}{119.6} & \textbf{-67.23} & \textcolor{revblue}{6.5} & \textbf{-67.23} & \textcolor{revblue}{\textbf{2.6}} \\
\bottomrule
\end{tabular}}
\caption{Results for \texttt{Baron} (time limit 3600 sec) , \texttt{Baron-H}, \textcolor{revblue}{\texttt{Knitro}, \texttt{Ipopt},} \AOsub and \texttt{LS$^*$} on polynomial kernels -  Breast Cancer Diagnostic dataset}
\label{table: baron results on bcd}
\end{table}

\begin{table}[!ht]
\centering
\renewcommand\arraystretch{1.2}\footnotesize
\resizebox{\textwidth}{!}{%
\begin{tabular}{ll|cc|cc|cc|cc|cc|cc}
\toprule
& & \multicolumn{2}{c|}{\texttt{Baron}} & \multicolumn{2}{c|}{\texttt{Baron-H}} & \multicolumn{2}{c|}{\textcolor{revblue}{\texttt{Knitro}}} & \multicolumn{2}{c|}{\textcolor{revblue}{\texttt{Ipopt}}} & \multicolumn{2}{c|}{\AOsub} & \multicolumn{2}{c}{\texttt{LS$^*$}}\\
 & & ObjFun & MipGap(\%)  & ObjFun & Time & \textcolor{revblue}{ObjFun} & \textcolor{revblue}{Time} & \textcolor{revblue}{ObjFun} & \textcolor{revblue}{Time} & ObjFun & \textcolor{revblue}{Time} & ObjFun & \textcolor{revblue}{Time} \\
\midrule
\multirow{5}{*}{\texttt{d=2}} & \texttt{m=10} & \textbf{-29.21} & 2.92e1 & -27.96 & \textbf{0.2} & \textcolor{revblue}{-23.68} & \textcolor{revblue}{1.0} & \textcolor{revblue}{-23.65} & \textcolor{revblue}{0.8} & -24.61 & \textcolor{revblue}{2.2} & \textbf{-29.21} & \textcolor{revblue}{2.3} \\
 & \texttt{m=50} & -86.74 & 1.28e3 & -86.74 & 5.5 & \textcolor{revblue}{-82.62} & \textcolor{revblue}{5.2} & \textcolor{revblue}{\textbf{-89.03}} & \textcolor{revblue}{5.2} & \textbf{-89.03} & \textcolor{revblue}{4.7} & \textbf{-89.03} & \textcolor{revblue}{\textbf{2.8}} \\
 & \texttt{m=100} & -108.56 & 4.94e3 & -108.56 & 59.6 & \textcolor{revblue}{-112.89} & \textcolor{revblue}{20.4} & \textcolor{revblue}{-112.89} & \textcolor{revblue}{22.7} & \textbf{-139.26} & \textcolor{revblue}{5.6} & \textbf{-139.26} & \textcolor{revblue}{\textbf{2.7}} \\
 & \texttt{m=150} & -178.70 & 6.86e4 & -178.70 & 246.1 & \textcolor{revblue}{-174.62} & \textcolor{revblue}{46.7} & \textcolor{revblue}{-172.22} & \textcolor{revblue}{75.4} & -190.70 & \textcolor{revblue}{4.2} & \textbf{-195.54} & \textcolor{revblue}{\textbf{3.8}} \\
 & \texttt{m=200} & -235.46 & 1.93e5 & -235.46 & 801.9 & \textcolor{revblue}{-229.57} & \textcolor{revblue}{87.8} & \textcolor{revblue}{-233.92} & \textcolor{revblue}{168.0} & \textbf{-244.64} & \textcolor{revblue}{6.3} & \textbf{-244.64} & \textcolor{revblue}{\textbf{4.0}} \\
\midrule
\multirow{5}{*}{\texttt{d=3}} & \texttt{m=10} & \textbf{-21.41} & 8.20e1 & -19.29 & \textbf{0.1} & \textcolor{revblue}{-19.29} & \textcolor{revblue}{0.4} & \textcolor{revblue}{-19.29} & \textcolor{revblue}{0.3} & -21.23 & \textcolor{revblue}{2.0} & -21.23 & \textcolor{revblue}{2.3} \\
 & \texttt{m=50} & \textbf{-52.54} & 6.13e3 & -49.46 & 6.2 & \textcolor{revblue}{-48.75} & \textcolor{revblue}{5.6} & \textcolor{revblue}{-48.75} & \textcolor{revblue}{6.1} & -50.89 & \textcolor{revblue}{4.2} & \textbf{-52.54} & \textcolor{revblue}{\textbf{2.9}} \\
 & \texttt{m=100} & -70.34 & 2.45e5 & -63.84 & 50.1 & \textcolor{revblue}{-74.05} & \textcolor{revblue}{20.5} & \textcolor{revblue}{-74.05} & \textcolor{revblue}{24.0} & \textbf{-75.32} & \textcolor{revblue}{3.2} & \textbf{-75.32} & \textcolor{revblue}{\textbf{2.9}} \\
 & \texttt{m=150} & -107.69 & 3.67e5 & -107.69 & 251.6 & \textcolor{revblue}{-103.65} & \textcolor{revblue}{46.8} & \textcolor{revblue}{-100.93} & \textcolor{revblue}{66.4} & -103.65 & \textcolor{revblue}{3.8} & \textbf{-108.66} & \textcolor{revblue}{\textbf{3.3}} \\
 & \texttt{m=200} & -136.19 & 1.01e6 & -136.19 & 996.4 & \textcolor{revblue}{-116.02} & \textcolor{revblue}{88.3} & \textcolor{revblue}{-115.88} & \textcolor{revblue}{162.1} & \textbf{-158.56} & \textcolor{revblue}{4.8} & \textbf{-158.56} & \textcolor{revblue}{\textbf{4.0}} \\
\midrule
\multirow{5}{*}{\texttt{d=5}} & \texttt{m=10} & \textbf{-18.12} & 4.46e2 & -6.68 & \textbf{0.1} & \textcolor{revblue}{\textbf{-18.12}} & \textcolor{revblue}{0.4} & \textcolor{revblue}{-7.45} & \textcolor{revblue}{0.4} & \textbf{-18.12} & \textcolor{revblue}{1.5} & \textbf{-18.12} & \textcolor{revblue}{2.4} \\
 & \texttt{m=50} & \textbf{-43.17} & 6.87e5 & \textbf{-43.17} & 5.3 & \textcolor{revblue}{-22.12} & \textcolor{revblue}{5.1} & \textcolor{revblue}{-21.79} & \textcolor{revblue}{5.2} & -29.84 & \textcolor{revblue}{3.0} & \textbf{-43.17} & \textcolor{revblue}{\textbf{2.8}} \\
 & \texttt{m=100} & \textbf{-45.93} & 4.02e6 & \textbf{-45.93} & 57.2 & \textcolor{revblue}{-21.01} & \textcolor{revblue}{19.9} & \textcolor{revblue}{-22.97} & \textcolor{revblue}{22.4} & -31.77 & \textcolor{revblue}{4.7} & -31.76 & \textcolor{revblue}{\textbf{3.0}} \\
 & \texttt{m=150} & -51.27 & 8.64e6 & -42.30 & 242.6 & \textcolor{revblue}{-35.38} & \textcolor{revblue}{47.4} & \textcolor{revblue}{-34.10} & \textcolor{revblue}{65.8} & -42.99 & \textcolor{revblue}{6.5} & \textbf{-71.79} & \textcolor{revblue}{\textbf{4.0}} \\
 & \texttt{m=200} & \textbf{-95.46} & 1.04e7 & \textbf{-95.46} & 808.0 & \textcolor{revblue}{-41.28} & \textcolor{revblue}{85.5} & \textcolor{revblue}{-42.12} & \textcolor{revblue}{155.6} & \textbf{-95.46} & \textcolor{revblue}{\textbf{4.1}} & \textbf{-95.46} & \textcolor{revblue}{4.4} \\
\bottomrule
\end{tabular}}
\caption{Results for \texttt{Baron} (time limit 3600 sec) , \texttt{Baron-H}, \textcolor{revblue}{\texttt{Knitro}, \texttt{Ipopt},} \AOsub and \texttt{LS$^*$} on polynomial kernels -  Ionosphere dataset}
\label{table: baron results iono}
\end{table}

\begin{table}[!ht]
\centering
\renewcommand\arraystretch{1.2}\footnotesize
\resizebox{\textwidth}{!}{%
\begin{tabular}{ll|cc|cc|cc|cc|cc|cc}
\toprule
& & \multicolumn{2}{c|}{\texttt{Baron}} & \multicolumn{2}{c|}{\texttt{Baron-H}} & \multicolumn{2}{c|}{\textcolor{revblue}{\texttt{Knitro}}} & \multicolumn{2}{c|}{\textcolor{revblue}{\texttt{Ipopt}}} & \multicolumn{2}{c|}{\AOsub} & \multicolumn{2}{c}{\texttt{LS$^*$}}\\
 & & ObjFun & MipGap(\%)  & ObjFun & Time & \textcolor{revblue}{ObjFun} & \textcolor{revblue}{Time} & \textcolor{revblue}{ObjFun} & \textcolor{revblue}{Time} & ObjFun & \textcolor{revblue}{Time} & ObjFun & \textcolor{revblue}{Time} \\
\midrule
\multirow{5}{*}{\texttt{d=2}} & \texttt{m=10} & \textbf{-2.35} & 2.93e3 & \textbf{-2.35} & \textbf{0.3} & \textcolor{revblue}{\textbf{-2.35}} & \textcolor{revblue}{0.5} & \textcolor{revblue}{\textbf{-2.35}} & \textcolor{revblue}{0.4} & \textbf{-2.35} & \textcolor{revblue}{3.9} & \textbf{-2.35} & \textcolor{revblue}{8.4} \\
 & \texttt{m=50} & -18.93 & 1.73e4 & -18.93 & 11.5 & \textcolor{revblue}{-18.46} & \textcolor{revblue}{9.7} & \textcolor{revblue}{-20.66} & \textcolor{revblue}{\textbf{8.7}} & -20.66 & \textcolor{revblue}{11.3} & \textbf{-20.78} & \textcolor{revblue}{9.6} \\
 & \texttt{m=100} & -39.23 & 3.27e5 & -39.23 & 118.5 & \textcolor{revblue}{-36.25} & \textcolor{revblue}{36.4} & \textcolor{revblue}{-43.50} & \textcolor{revblue}{39.2} & \textbf{-47.11} & \textcolor{revblue}{\textbf{11.8}} & \textbf{-47.11} & \textcolor{revblue}{14.6} \\
 & \texttt{m=150} & -35.69 & 7.87e5 & -35.69 & 453.5 & \textcolor{revblue}{-55.29} & \textcolor{revblue}{88.8} & \textcolor{revblue}{-58.53} & \textcolor{revblue}{106.6} & \textbf{-69.20} & \textcolor{revblue}{19.7} & \textbf{-69.20} & \textcolor{revblue}{\textbf{17.3}} \\
 & \texttt{m=200} & -72.79 & 4.68e5 & -72.79 & 696.9 & \textcolor{revblue}{-53.04} & \textcolor{revblue}{101.1} & \textcolor{revblue}{-56.20} & \textcolor{revblue}{128.6} & \textbf{-78.49} & \textcolor{revblue}{23.9} & \textbf{-78.49} & \textcolor{revblue}{\textbf{23.0}} \\
\midrule
\multirow{5}{*}{\texttt{d=3}} & \texttt{m=10} & \textbf{-0.82} & 4.24e4 & -0.81 & \textbf{0.2} & \textcolor{revblue}{\textbf{-0.82}} & \textcolor{revblue}{0.5} & \textcolor{revblue}{\textbf{-0.82}} & \textcolor{revblue}{1.1} & \textbf{-0.82} & \textcolor{revblue}{3.6} & \textbf{-0.82} & \textcolor{revblue}{9.9} \\
 & \texttt{m=50} & -3.94 & 4.32e6 & -3.94 & 12.2 & \textcolor{revblue}{-3.92} & \textcolor{revblue}{9.5} & \textcolor{revblue}{-4.19} & \textcolor{revblue}{8.6} & \textbf{-4.26} & \textcolor{revblue}{\textbf{6.6}} & \textbf{-4.26} & \textcolor{revblue}{10.6} \\
 & \texttt{m=100} & -6.79 & 8.10e6 & -6.79 & 97.3 & \textcolor{revblue}{-7.60} & \textcolor{revblue}{35.1} & \textcolor{revblue}{-8.08} & \textcolor{revblue}{36.6} & \textbf{-8.93} & \textcolor{revblue}{21.7} & \textbf{-8.93} & \textcolor{revblue}{\textbf{14.7}} \\
 & \texttt{m=150} & -9.49 & 1.26e7 & -9.49 & 496.0 & \textcolor{revblue}{-10.81} & \textcolor{revblue}{88.9} & \textcolor{revblue}{-10.60} & \textcolor{revblue}{90.8} & \textbf{-11.18} & \textcolor{revblue}{\textbf{12.5}} & \textbf{-11.18} & \textcolor{revblue}{19.9} \\
 & \texttt{m=200} & -9.96 & 1.42e7 & -9.96 & 687.8 & \textcolor{revblue}{-11.30} & \textcolor{revblue}{103.7} & \textcolor{revblue}{-11.51} & \textcolor{revblue}{117.5} & \textbf{-12.37} & \textcolor{revblue}{\textbf{15.7}} & \textbf{-12.37} & \textcolor{revblue}{22.1} \\
\midrule
\multirow{5}{*}{\texttt{d=5}} & \texttt{m=10} & \textbf{-0.20} & 8.10e7 & \textbf{-0.20} & \textbf{0.2} & \textcolor{revblue}{\textbf{-0.20}} & \textcolor{revblue}{0.7} & \textcolor{revblue}{\textbf{-0.20}} & \textcolor{revblue}{0.4} & \textbf{-0.20} & \textcolor{revblue}{2.1} & \textbf{-0.20} & \textcolor{revblue}{8.4} \\
 & \texttt{m=50} & -0.63 & 1.92e9 &-0.62 & 12.0 & \textcolor{revblue}{\textbf{-0.66}} & \textcolor{revblue}{10.0} & \textcolor{revblue}{-0.56} & \textcolor{revblue}{9.7} &  \textbf{-0.66} & \textcolor{revblue}{\textbf{8.8}} &  \textbf{-0.66} & \textcolor{revblue}{10.8} \\
 & \texttt{m=100} & -0.99 & 2.78e9 & -0.99 & 119.3 & \textcolor{revblue}{-0.99} & \textcolor{revblue}{36.9} & \textcolor{revblue}{-0.99} & \textcolor{revblue}{40.1} & -1.02 & \textcolor{revblue}{17.2} & \textbf{-1.08} & \textcolor{revblue}{\textbf{16.6}} \\
 & \texttt{m=150} & -1.21 & 4.03e9 & -1.21 & 490.6 & \textcolor{revblue}{-1.19} & \textcolor{revblue}{88.4} & \textcolor{revblue}{-1.22} & \textcolor{revblue}{105.7} & \textbf{-1.24} & \textcolor{revblue}{\textbf{13.7}} & \textbf{-1.24} & \textcolor{revblue}{23.3} \\
 & \texttt{m=200} & -1.30 & 6.27e9 & -1.30 & 747.7 & \textcolor{revblue}{-1.32} & \textcolor{revblue}{122.8} & \textcolor{revblue}{-1.30} & \textcolor{revblue}{134.9} & \textbf{-1.33} & \textcolor{revblue}{\textbf{14.9}} & \textbf{-1.33} & \textcolor{revblue}{24.9} \\
\bottomrule
\end{tabular}}
\caption{Results for \texttt{Baron} (time limit 3600 sec) , \texttt{Baron-H}, \textcolor{revblue}{\texttt{Knitro}, \texttt{Ipopt},} \AOsub and \texttt{LS$^*$} on polynomial kernels - Sonar dataset}
\label{table: baron results sonar}
\end{table}

\begin{table}[!ht]
\centering
\renewcommand\arraystretch{1.2}\footnotesize
\resizebox{\textwidth}{!}{%
\begin{tabular}{ll|cc|cc|cc|cc|cc}
\toprule
& & \multicolumn{2}{c|}{\texttt{Baron}} & \multicolumn{2}{c|}{\texttt{Baron-H}} & \multicolumn{2}{c|}{\textcolor{revblue}{\texttt{Knitro}}} & \multicolumn{2}{c|}{\textcolor{revblue}{\texttt{Ipopt}}} & \multicolumn{2}{c}{\texttt{LS$^*$}}\\
 & & ObjFun & MipGap(\%)  & ObjFun & Time & \textcolor{revblue}{ObjFun} & \textcolor{revblue}{Time} & \textcolor{revblue}{ObjFun} & \textcolor{revblue}{Time} & ObjFun & \textcolor{revblue}{Time}  \\
\midrule
\multirow{5}{*}{Breast Cancer D.} & \texttt{m=10} & \textbf{-4.70} & 2.51e4 & \textbf{-4.70} & \textbf{0.05} & \textcolor{revblue}{\textbf{-4.70}} & \textcolor{revblue}{0.30} & \textcolor{revblue}{\textbf{-4.70}} & \textcolor{revblue}{1.00} & \textbf{-4.70} & \textcolor{revblue}{1.93} \\
 & \texttt{m=50} & \textbf{-55.32} & 5.14e5 & \textbf{-55.32} & 3.12 & \textcolor{revblue}{-48.71} & \textcolor{revblue}{4.70} & \textcolor{revblue}{-48.71} & \textcolor{revblue}{3.20} & \textbf{-55.32} & \textcolor{revblue}{\textbf{2.37}} \\
 & \texttt{m=100} & -96.87 & 1.15e5 & -96.87 & 41.13 & \textcolor{revblue}{-98.42} & \textcolor{revblue}{16.80} & \textcolor{revblue}{-98.42} & \textcolor{revblue}{13.70} & \textbf{-98.82} & \textcolor{revblue}{\textbf{2.58}} \\
 & \texttt{m=150} & -114.27 & 2.36e6 & -114.27 & 199.44 & \textcolor{revblue}{-120.12} & \textcolor{revblue}{39.00} & \textcolor{revblue}{-117.48} & \textcolor{revblue}{33.10} & \textbf{-141.04} & \textcolor{revblue}{\textbf{3.10}} \\
 & \texttt{m=200} & \textbf{-190.70} & 4.92e6 & \textbf{-190.70} & 706.96 & \textcolor{revblue}{-161.45} & \textcolor{revblue}{69.30} & \textcolor{revblue}{-161.45} & \textcolor{revblue}{68.00} & \textbf{-190.70} & \textcolor{revblue}{\textbf{3.66}} \\
\midrule
\multirow{5}{*}{Cleveland} & \texttt{m=10} & \textbf{-39.90} & 2.50e3 & \textbf{-39.90} & \textbf{0.06} & \textcolor{revblue}{\textbf{-39.90}} & \textcolor{revblue}{0.20} & \textcolor{revblue}{\textbf{-39.90}} & \textcolor{revblue}{0.20} & \textbf{-39.90} & \textcolor{revblue}{0.98} \\
 & \texttt{m=50} & -215.46 & 1.19e5 & -207.66 & \textbf{0.79} & \textcolor{revblue}{-195.87} & \textcolor{revblue}{2.30} & \textcolor{revblue}{-207.66} & \textcolor{revblue}{2.20} & \textbf{-215.46} & \textcolor{revblue}{1.08} \\
 & \texttt{m=100} & -487.66 & 3.60e5 & -487.66 & 9.94 & \textcolor{revblue}{-483.95} & \textcolor{revblue}{7.60} & \textcolor{revblue}{-487.66} & \textcolor{revblue}{6.90} & \textbf{-520.96} & \textcolor{revblue}{\textbf{1.65}} \\
 & \texttt{m=150} & \textbf{-780.49} & 3.19e5 & \textbf{-780.49} & 47.32 & \textcolor{revblue}{-743.03} & \textcolor{revblue}{18.30} & \textcolor{revblue}{-712.86} & \textcolor{revblue}{17.50} & \textbf{-780.49} & \textcolor{revblue}{\textbf{1.88}} \\
 & \texttt{m=200} & -1106.81 & 9.00e5 & -1106.81 & 151.27 & \textcolor{revblue}{-1085.60} & \textcolor{revblue}{34.70} & \textcolor{revblue}{-1106.81} & \textcolor{revblue}{34.10} & \textbf{-1113.48} & \textcolor{revblue}{\textbf{2.90}} \\
\midrule
\multirow{5}{*}{Ionosphere} & \texttt{m=10} & \textbf{-32.20} & 5.48e2 & \textbf{-32.20} & \textbf{0.06} & \textcolor{revblue}{-25.79} & \textcolor{revblue}{0.40} & \textcolor{revblue}{-25.55} & \textcolor{revblue}{0.30} & \textbf{-32.20} & \textcolor{revblue}{2.48} \\
 & \texttt{m=50} & -85.42 & 3.25e5 & -85.42 & \textbf{2.96} & \textcolor{revblue}{-84.22} & \textcolor{revblue}{5.20} & \textcolor{revblue}{-84.22} & \textcolor{revblue}{4.00} & \textbf{-87.36} & \textcolor{revblue}{3.27} \\
 & \texttt{m=100} & -124.26 & 9.10e5 & -124.26 & 40.25 & \textcolor{revblue}{-113.95} & \textcolor{revblue}{18.00} & \textcolor{revblue}{-113.95} & \textcolor{revblue}{15.00} & \textbf{-134.43} & \textcolor{revblue}{\textbf{3.67}} \\
 & \texttt{m=150} & -155.69 & 1.73e6 & -155.69 & 195.34 & \textcolor{revblue}{-147.60} & \textcolor{revblue}{41.20} & \textcolor{revblue}{-147.15} & \textcolor{revblue}{39.30} & \textbf{-183.86} & \textcolor{revblue}{\textbf{4.62}} \\
 & \texttt{m=200} & -209.81 & 4.52e6 & -209.81 & 658.14 & \textcolor{revblue}{-176.07} & \textcolor{revblue}{76.90} & \textcolor{revblue}{-177.94} & \textcolor{revblue}{78.10} & \textbf{-216.25} & \textcolor{revblue}{\textbf{6.09}} \\
\midrule
\multirow{5}{*}{Sonar} & \texttt{m=10} & \textbf{-6.80} & 4.34e4 & -6.01 & \textbf{0.12} & \textcolor{revblue}{-6.01} & \textcolor{revblue}{0.50} & \textcolor{revblue}{\textbf{-6.80}} & \textcolor{revblue}{0.40} & \textbf{-6.80} & \textcolor{revblue}{8.47} \\
 & \texttt{m=50} & -39.31 & 8.16e5 & -39.24 & 9.93 & \textcolor{revblue}{-36.75} & \textcolor{revblue}{8.10} & \textcolor{revblue}{-36.82} & \textcolor{revblue}{\textbf{6.40}} & \textbf{-39.32} & \textcolor{revblue}{12.54} \\
 & \texttt{m=100} & -65.59 & 2.07e6 & -65.59 & 84.25 & \textcolor{revblue}{-64.65} & \textcolor{revblue}{29.80} & \textcolor{revblue}{-64.61} & \textcolor{revblue}{25.00} & \textbf{-66.31} & \textcolor{revblue}{\textbf{21.75}} \\
 & \texttt{m=150} & -87.75 & 3.44e6 & -87.75 & 405.96 & \textcolor{revblue}{-90.13} & \textcolor{revblue}{68.30} & \textcolor{revblue}{-90.13} & \textcolor{revblue}{63.60} & \textbf{-91.18} & \textcolor{revblue}{\textbf{36.66}} \\
 & \texttt{m=200} & -96.94 & 7.08e6 & -96.94 & 623.11 & \textcolor{revblue}{-99.68} & \textcolor{revblue}{85.30} & \textcolor{revblue}{-99.68} & \textcolor{revblue}{74.10} & \textbf{-101.85} & \textcolor{revblue}{\textbf{41.39}} \\
\bottomrule
\end{tabular}}
\caption{Results for \texttt{Baron} (time limit 3600 sec), \texttt{Baron-H}, \textcolor{revblue}{\texttt{Knitro}, \texttt{Ipopt},} and \texttt{LS*} across various datasets on a Gaussian kernel}
\label{table: baron results gaussian kernel}
\end{table}

Tables~\ref{table: baron results clev}, \ref{table: baron results on bcd}, \ref{table: baron results iono}, and~\ref{table: baron results sonar} report the results obtained using different polynomial kernels, while Table~\ref{table: baron results gaussian kernel} refers to the case where a Gaussian kernel is used.

In general, the computational time required by \texttt{Baron-H} increases sharply as the number of data points $m$ grows. Regarding the standard resolution process of \texttt{Baron}, a time limit of 3600 seconds was imposed, but this limit was always reached, except for one case with the Cleveland dataset for $d=5$ and $m=10$. Moreover, the MIP gap values are generally surprisingly high, and tend to increase further with the degree of the polynomial kernel. Although the computational times of our algorithms also increase with $m$, the growth is much less drastic compared to \texttt{Baron-H}. In addition, both \AOsub and \texttt{LS$^*$} almost always find the best solutions in the shortest computational times, particularly for higher, and more realistic values of $m$.

\end{document}